\documentclass{article}
\usepackage{graphicx} % Required for inserting images
\usepackage{amssymb}
\usepackage{yhmath}
\usepackage{amsmath}
\usepackage{amsthm}
\usepackage{caption}
\usepackage{subcaption}
\usepackage[left=1 in,right=1 in,top=1 in,bottom=1 in]{geometry}
\usepackage[round]{natbib}
\usepackage[colorlinks=true]{hyperref}
\usepackage{xfrac}
\newtheorem{thm}{Theorem}
\newtheorem{crl}{Corollary}
\newtheorem{lmm}{Lemma}

\allowdisplaybreaks

\usepackage{bigfoot}

\DeclareNewFootnote{AAffil}[arabic]
\DeclareNewFootnote{ANote}[fnsymbol]

\usepackage{etoolbox}
\makeatletter
\patchcmd\maketitle{\def\@makefnmark{\rlap{\@textsuperscript{\normalfont\@thefnmark}}}}{}{}{}
\makeatother

\makeatletter
\def\thanksAAffil#1{% <--- These %'s are necessary for spacing
  \footnotemarkAAffil\protected@xdef\@thanks{\@thanks%
        \protect\footnotetextAAffil[\the \c@footnoteAAffil]{#1}}%
}
\def\thanksANote#1{%
  \footnotemarkANote%
  \protected@xdef\@thanks{\@thanks%
        \protect\footnotetextANote[\the \c@footnoteANote]{#1}}%
}
\makeatother
\begin{document}
\title{On Vector Field Reconstruction from Noisy ODE in High\\ Ambient Dimension}
\author{Hugo Henneuse \thanksAAffil{Centre Borelli, ENS Paris Saclay, Université Paris-Saclay, France}$^{\text{ ,}}$\thanksAAffil{DataShape, Inria Saclay, Palaiseau, France}\\ \href{mailto:hugo.henneuse@universite-paris-saclay.fr}{hugo.henneuse@universite-paris-saclay.fr}}
\maketitle
\begin{abstract}
This work investigates the nonparametric estimation of the vector field of a noisy Ordinary Differential Equation (ODE) in high-dimensional ambient spaces, under the assumption that the initial conditions are sampled from a lower-dimensional structure. Specifically, let \( f:\mathbb{R}^{D}\to\mathbb{R}^{D} \) denote the vector field of the autonomous ODE \( y' = f(y) \). We observe noisy trajectories \( \tilde{y}_{X_i}(t_j) = y_{X_i}(t_j) + \varepsilon_{i,j} \), where \( y_{X_i}(t_j) \) is the solution at time \( t_j \) with initial condition \( y(0)=X_i \), the \( X_i \) are drawn from a \((a,b)\)-standard distribution \( \mu \), and \( \varepsilon_{i,j} \) denotes noise. From a minimax perspective, we study the reconstruction of \( f \) within the envelope of trajectories generated by the support of \( \mu \). We proposed an estimator combining flow reconstruction with derivative estimation techniques from nonparametric regression. Under mild regularity assumptions on \( f \), we establish convergence rates that depend on the temporal resolution, the number of initial conditions, and the parameter \( b \), which controls the mass concentration of \( \mu \). These rates are then shown to be minimax optimal (up to logarithmic factors) and illustrate how the proposed approach mitigates the curse of dimensionality. Additionally, we illustrate the computational and statistical efficiency of our estimator through numerical experiments.
\end{abstract}
\section*{Introduction}
Longitudinal data commonly arise in the study of dynamical systems across various disciplines, from physics to biology and economics \citep[see e.g.][]{Hirsch2013}. By observing multiple trajectories in a given space, we aim to extract key information about the underlying dynamics, such as governing laws, stability properties, or long-term trends. A common approach to modeling such phenomena is to assume that the observed trajectories follow an underlying differential system, which provides a structured framework for inference, prediction, and control.

However, real-world data are often subject to noise, missing observations, and irregular sampling, making their analysis particularly challenging and raising interesting statistical questions. In particular, given noisy observations of trajectories governed by an autonomous differential equation:
\begin{equation}
\label{eq:ODE}
    y' = f(y),
\end{equation}
where \( f: \mathbb{R}^D \to \mathbb{R}^D \) is the vector field governing the system, can we reconstruct \( f \) from the observed data?

This question has garnered significant interest among the machine learning and statistical communities. Most existing works focus on the case where $f$ is assumed to belong to a parametric class of functions \citep[e.g., see the surveys][]{McgoffSurvey15, RamsaySurvey17, DattnerSurvey21}. In contrast, the non-parametric literature is more limited and primarily numerical. Historically, popular methods are based on ``gradient matching'' techniques \citep{Brunton2016, Sneed2021}, which first estimate the gradient from the data and then perform nonparametric regression to recover $f$. Several later approaches model $f$ using deep learning \citep{Chen18, RUDY2019, Bhat20, GOTTWALD2021, Marzouk2024}. Alternatively, \cite{Heinonen18} proposes a Bayesian framework using Gaussian processes to model the vector field. More recently, \cite{LAHOUEL2024} introduced an estimator based on solving a constrained optimization problem in an RKHS. Some of these papers provide theoretical guarantees, but mostly in a density setting \citep[for instance, see][Theorem 2.2]{Marzouk2024}, which is distinct from the regression framework we consider. Closer to our setting, \cite{LAHOUEL2024} provides consistency guarantees (see Corollary 1), but only for trajectory estimation. The overall lack of statistical guarantees for the estimation of the vector field in regression settings motivates a theoretical study of this problem.

In this direction, the two recent papers \cite{schötz2024, schötz2025} provide important theoretical insights into the problem. More precisely, they introduce two statistical models: the \emph{snake model} and the \emph{stubble model}. In the snake model, a few long trajectories are observed, covering the space of interest. In the stubble model, many short trajectories are observed, with initial values forming a deterministic and uniform cover of the space of interest. In both cases, the authors propose regression-based estimation strategies and provide a detailed analysis of their convergence rates, depending on the regularity of \( f \). To our knowledge, these works represent the most advanced contributions to the statistical understanding of the vector field reconstruction problem.\\\\
In this paper, we propose a different framework that we believe is better suited for high-dimensional applications. Our model can be interpreted as a simple regression model on the flow, with a random design over the initial values. More precisely, we assume that the initial values are sampled from a probability distribution \( \mu \) with compact support \( \mathcal{X} \) and that we observe noisy trajectories from these sampled initial conditions. In this setting, we address the problem of inferring \( f \) within the envelope of solutions originating from \( \mathcal{X} \), which defines our space of interest. Unlike previous models, particularly those considered in \cite{schötz2024, schötz2025}, our framework introduces significant flexibility to handle such scenarios. The observed initial values are not required to cover the space of interest. Moreover, they are random and may produce a highly non-uniform covering of \( \mathcal{X} \), as we will consider very general classes of probability distributions. Consequently, the observed trajectories do not necessarily provide a uniform covering of the space of interest.

A key motivation for considering this model arises from experimental settings in which the evolutions of individuals are recorded in a parameter space of potentially high dimension, while the initial positions (and, more generally, the resulting trajectories) are confined to a much lower-dimensional structure (for instance a sub-manifold of the ambient space), often as a consequence of (possibly unknown) constraints. This perspective has recently attracted growing interest. The main line of work \citep{Linot2020, Fukami2021, Vlachas2022, Wu2024, Regazzoni2024} seeks to reduce dimensionality by identifying a lower-dimensional latent space that captures the global dynamics, and then learning the dynamics within this latent space. A key limitation of these approaches is that the learned latent space often has a significantly higher dimension than the intrinsic dimension of the trajectories, and its construction is computationally expensive. An alternative was recently proposed by \cite{Huang2025}, who adapt the methods of \cite{Brunton2016} to reconstruct the vector field of an ODE under the geometric constraint that trajectories evolve on a submanifold of the ambient space. Their approach estimates local tangent spaces of the underlying manifold and uses them to construct a basis of functions for approximating the vector field. Unlike the aforementioned latent-space methods, this construction ensures that the estimated vector field remains close to the tangent spaces of the manifold, so that the resulting trajectories do not deviate substantially from it. However, these works are primarily numerical: although some results are provided on the consistency of the estimated trajectories \citep[see, e.g., Theorem~3.1 in][]{Huang2025}, no theoretical guarantees have yet been established for vector field estimation itself.\\\\
 \textbf{Contributions.} In this context, our contributions are as follows :
 \begin{itemize}
     \item We study the statistical complexity of the problem from a minimax perspective by establishing lower bounds on pointwise estimation rates. This analysis encompasses broad classes of initial value distributions and operates under mild assumptions about the vector field, revealing three distinct regimes based on temporal resolution, the number of observed trajectories, and the mass concentration of the initial values distribution. Importantly, those rates do not depend on the ambient dimension.
     \item We propose a simple estimation procedure, based on derivative estimation techniques from nonparametric regression combined with nearest-neighbor averaging, which attains the obtained lower bounds (up to logarithmic factors) and thereby establishes their optimality. Moreover, although our method may not outperform more sophisticated estimators in practice, it enjoys several desirable properties: it has low computational complexity, requires minimal prior knowledge about the structure supporting the sampled trajectories, and does not rely on explicitly estimating this structure.
     \item We explore connections with manifold learning to clarify how our approach mitigates the curse of dimensionality, especially when trajectories evolve in high-dimensional spaces with initial values lying on a low-dimensional submanifold.
     \item We additionally illustrate that our analysis naturally extends to other common risk metrics, establishing that our estimation procedure achieves (up to logarithmic factors) the same rates for the sup-norm risk.
 \end{itemize}
\textbf{Organization.} The paper is organized as follows. Section \ref{sec: frmwk,mthd,rslt} specifies the formal framework of this work and describes our estimation strategy. Section \ref{sec: main results} states our general theoretical results, and Section \ref{sec: geom-coro} derives their specific forms under manifold assumptions. Section \ref{sec: numerical illustration} presents numerical illustrations of our results and examines the computational aspects of the proposed estimator. Detailed proofs of all results, along with auxiliary lemmas and an extension to convergence in the sup-norm, are deferred to the Appendix.
\section{Framework and method}
\label{sec: frmwk,mthd,rslt}

We now present the precise framework and estimation strategy considered in this work, followed by our main convergence results. 

\subsection{Framework}

\textbf{Assumption on the vector field.} We consider the autonomous Ordinary Differential Equation (ODE) defined by (\ref{eq:ODE}) and assume that \( f \) belongs to \( \operatorname{Lip}(L,M) \), the set of \( L \)-Lipschitz functions from \( \mathbb{R}^{D} \) to \( \mathbb{R}^{D} \) with norm bounded by \( M \). Let \( \mathcal{X} \subset \mathbb{R}^{D} \) be the space of initial values. By the Cauchy-Lipschitz theorem, for any \( x \in \mathbb{R}^{D} \), there exists a unique global solution \( y_x \) to (\ref{eq:ODE}) with initial condition \( y(0) = x \). For all \( t \in [0,+\infty[ \), we denote by \( \phi: (x,t) \mapsto y_x(t) \) the flow associated with (\ref{eq:ODE}). Since \( f \) is Lipschitz, \( \phi \) is continuously differentiable.\\\\
\textbf{Assumption on the initial values sampling.} Let \( \mu \) be a probability distribution with support \( \mathcal{X} \) from which will be sampled the initial values. We suppose that \(\mu\) belongs, for some \(0<a, 0 \leq b\leq D-1\), to \( \mathcal{P}(a,b) \), the class of \( (a,b) \)-standard probability distributions, i.e., the set of probability distributions \( \mu \) such that for all \( x \) in its support \(\operatorname{supp}(\mu) \) and \( R > 0 \), we have:
\[
\mu(B_2(x,R)) \geq \min(a R^b, 1),
\]
where \( B_2(x,R) \) denotes the Euclidean ball of radius \( R \) centered at \( x \). This condition arises originally from the (level) set estimation literature \citep[see e.g.][]{Cuevas2004,Cuevas2009,ChazalGlisseMichel} and encompasses large classes of distributions. It gives, at the same time, lower control on the mass of $\mu$ and (in some sense) on the geometry of \(\operatorname{supp}(\mu) \). In some sense, the parameter $b$ can be linked to the notion of ``intrinsic dimension'' of $\mathcal{X}$. This connection becomes particularly clear when the initial conditions are sampled from a finite-dimensional manifold (see Section~\ref{sec: geom-coro}). However, it is worth emphasizing that the $(a,b)$-standard assumption imposes very little regularity on the support: it need not be a manifold and may even be highly singular, allowing for fractal-like structures.\\\\
\textbf{Regression settings.} We observe \( X = \{X_1, \dots, X_n\} \subset \mathcal{X} \), a set of \( n \) points sampled independently from \( \mu \) and we consider $G_T$ a regular temporal grid \( 0 < t_1 < t_2 < \dots < t_m = T \) over the interval \( [0,T] \). The observed discrete and noisy trajectories are given by:
\[
\tilde{y}_{X_i}(t_j) = y_{X_i}(t_j) + \sigma \varepsilon_{i,j},
\]
where \( \sigma > 0 \) and \( (\varepsilon_{i,j})_{i \in \{1,...,n\}, j \in \{1,...,m\}} \) are independent standard \( D \)-multivariate Gaussian variables, assumed to be independent of \( X \). Reformulating in the nonparametric regression framework, we observe:
\[
Y_{i,j} = \tilde{\phi}(X_i, t_j) = \phi(X_i, t_j) + \sigma \varepsilon_{i,j}, \quad \forall i \in \{1,...,n\}, j \in \{1,...,m\}.
\]
\textbf{Goal.} Our objective is to recover \( f(x) \) for any \( x \in \operatorname{Env}_{f}(\mathcal{X},T) \), where:
\[
\operatorname{Env}_{f}(\mathcal{X},T) := \{\phi(x,t) \mid x \in \mathcal{X}, t \in [0,T]\}.
\]
From a regression perspective, since there exist \( \tilde{x} \in \mathcal{X} \) and \( t \in [0,T] \) such that \( \phi(\tilde{x},t) = x \), this problem is equivalent to estimating the partial derivative \( \frac{d}{dt} \phi(\tilde{x},t) \). However, as we do not have direct access to the latent parametrization (\(\tilde{x}, t\)), standard regression-based derivative estimation techniques cannot be directly applied. This represents the main challenge when dealing with ``ODE-regression'' models.

\subsection{Estimator}

We now propose a three-steps estimator for \( f \). First, split the data in half spatially, i.e. consider $X[1,\lfloor n/2\rfloor]=\{X_{1},...,X_{\lfloor n/2\rfloor}\}$ and $X[\lfloor n/2\rfloor,n]=\{X_{\lfloor n/2\rfloor+1},...,X_{n}\}$ and then :
\begin{itemize}
    \item \textbf{Step 1 (Flow Estimation):} we first construct an estimator of the flow \( \phi \) using spatio-temporal nearest neighbors from the trajectories starting at points from  $X[1,\lfloor n/2\rfloor]$:
    \[
    \hat{\phi}(z,t) = \frac{1}{k_1 k_2} \sum_{p=1}^{k_1} \sum_{q=1}^{k_2} \tilde{\phi}(X^p(z), t^q(t)),
    \]
    where \( X^1(z), \dots, X^{k_1}(z) \) are the \( k_1 \)-nearest neighbors of \( z \) in \( X[1,\lfloor n/2\rfloor] \), and \( t^1(t), \dots, t^{k_2}(t) \) are the (deterministic) \( k_2 \)-nearest neighbors of \( t \) in $\{t_{1},...,t_{m}\}$. Note that this choice of flow estimator is somewhat arbitrary and can be replaced by any regression estimator satisfying Lemma \ref{lmm: regression flow} (or a similar condition). The key property we will use in this work is that the selected flow estimator achieves faster rates than those stated in Theorem \ref{thm: estim-lip-ab}.
    
    \item \textbf{Step 2 (Derivative Estimation):} we then estimate the vector field \( f \) along the trajectories starting at points in \(Z= X[\lfloor n/2\rfloor,n] \), adapting the strategy from \cite{Brabanter13} for derivative estimation in the non-parametric regression setting. We use the following derivative estimators :
    \[
    \widehat{f(\phi(X_i,t_j))} = \sum_{l=1}^{k} w_l \left(\frac{\tilde{\phi}(X_i,t_{j+l}) - \tilde{\phi}(X_i,t_{j-l})}{t_{j+l} - t_{j-l}}\right),
    \]
    for all \( X_i \in Z \) and \( j \in [k, m-k] \), where \( k > 0 \) is a parameter and the weights \( w_l \) are defined as:
    \[
    w_l = \frac{6l^2}{k(k+1)(2k+1)}.
    \]
    This choice of weights minimizes the variance of the derivative estimator under the constraint that their sum equals 1. Further details can be found in the proof of Theorem 1 in \cite{Brabanter13}. We omit derivative estimation near the temporal boundaries \( t = 0 \) and \( t = T \) as the loss of information is negligible for achieving optimal convergence rates. 
    
    \item \textbf{Step 3 (Vector Field Reconstruction) :} finally, we derive an estimator of \( f(x) \) for \( x \in \operatorname{Env}_{f}(\mathcal{X},T) \) using nearest-neighbor averaging :
    \[
    \hat{f}(x) = \frac{1}{r} \sum_{l=1}^{r} \widehat{f(\phi(X_{i_{l}(x)}, t_{j_{l}}(x)))},
    \]
    where :
    \[
    X_{i_{1}(x)} \in \operatorname{argmin} \left\{\min_{k+1 \leq j \leq m-k-1, t_{j}\in G_{T}} \|\hat{\phi}(X_i, t_j) - x\|, X_i \in Z\right\},
    \]
    and for \( 2 \leq l \leq r \),
    \[
    X_{i_{l}(x)} \in \operatorname{argmin} \left\{\min_{k+1 \leq j \leq m-k-1, t_{j}\in G_{T}} \|\hat{\phi}(X_i, t_j) - x\|, X_i \in Z \setminus \{X_{i_{1}(x)}, \dots, X_{i_{l-1}(x)}\} \right\}.
    \]
    The corresponding times \( t_{j_{l}(x)} \) are chosen as:
    \[
     t_{j_{l}(x)} \in \operatorname{argmin}_{k+1 \leq j \leq m-k-1,t_{j}\in G_{T}} \|x - \phi(X_{i_{l}(x)}, t_j)\|.
    \]
More informally, this final step consists of determining $(\hat{\phi}(X_{i_{l}(x)}, t_{j_{l}}(x))))_{1\leq l\leq r}$ the $r-$closest points from estimated trajectories among the ones starting on $Z$ and averaging the estimated derivatives of the associated trajectories at time $t_{j_{l}}$.
\end{itemize}
\textbf{Remark 1.} Observe that our estimator relies on simple tools: weighted finite differences for Step~2 and $k$-nearest-neighbors regression for Steps~1 and~3. Although it is not intended to outperform more sophisticated methods, such as those cited in the introduction, its simplicity makes it particularly suitable for theoretical analysis, allowing a precise quantification of convergence rates. This simplicity also yields relatively low computational complexity: linear in the number of observations $nm$ and the ambient dimension $D$.  

For practical use, several refinements could enhance performance. For example, to address boundary issues in Step~2, one could adapt strategies from \cite{Brabanter13}, such as those discussed in Section~2.2. When dealing with vector fields that are smoother than Lipschitz-continuous, weighted finite differences can be replaced by smoother variants, as proposed in \cite{DAI2016}. Similarly, for Steps~1 and~3, local averaging could be replaced by higher-order methods, such as local polynomial or kernel regression. While these modifications may improve statistical performance, they can also substantially increase computational cost, highlighting an inherent trade-off.

\section{Main Results}\label{sec: main results}
Under our framework, we show that for appropriately calibrated choices of $k_{1}$, $k_{2}$, \( k \) and \( r \), the proposed estimator achieves the following convergence rate.
\begin{thm}
\label{thm: estim-lip-ab}
For the choices of calibration parameters :
\begin{itemize}
    \item $k_{1}\simeq \max\left(\frac{(n/\ln(n))^{3/(b+3)}}{(m/\ln(nm))^{b/(b+3)}},1\right)$
    \item $k_{2}\simeq\max\left( \frac{m^{(b+2)/(b+3)}}{(n/(\ln(n)\ln(nm)))^{1/(b+3)}},1\right)$
    \item  \( r \simeq \max \left( \frac{n^{5/(5+b)}}{m^{b/(5+b)}}, 1 \right) \)
    \item \( k \simeq \max \left( \frac{m^{(4+b)/(5+b)}}{n^{1/(5+b)}}, 1 \right) \)
\end{itemize}
and $n$ and $m$ sufficiently large, we have:
$$\sup_{\substack{f \in \operatorname{Lip}(L,M)\\\mu \in \mathcal{P}(a,b)}} \sup_{x \in \operatorname{Env}_f(\operatorname{supp}(\mu), T)} \mathbb{E} \left[ \|\hat{f}(x) - f(x)\| \right]\leq C_{1} \max \left( \left( \frac{1}{nm} \right)^{\frac{1}{5+b}}, \left( \frac{\ln(n)}{n} \right)^{\frac{1}{b}}, \frac{1}{m} \right),$$
where \( C_{1} \) is a constant depending only on \( L \), $M$, \( a \), \( b \), \( D \), \( T \), and \( \sigma \).
\end{thm}

First, note that the convergence rates we obtain do not depend on the ambient dimension $D$, but rather on the parameter $\beta$, which controls the mass concentration of $\mu$ and can be, in many cases, linked to the ``dimension'' of its support (see Section \ref{sec: geom-coro}). This is a particularly interesting result in the context of dynamical systems evolving in high-dimensional ambient spaces, when the initial conditions (and hence the trajectories) lie on a lower-dimensional structure. Second, the result also highlights three distinct convergence regimes, depending on the values of $m$ and $n$. In light of the proof of Theorem~\ref{thm: estim-lip-ab} (see Section~\ref{sec:proof 1+crl}), these regimes can be interpreted as follows:
\begin{itemize}
    \item a ``trajectory-rich'' regime when :
    $$\max \left( \left( \frac{1}{nm} \right)^{\frac{1}{5+b}}, \left( \frac{\ln(n)}{n} \right)^{\frac{1}{b}}, \frac{1}{m} \right)=\frac{1}{m}.$$ 
    In this case, the convergence rate is \( O\left( \frac{1}{m} \right) \), which essentially corresponds to the bias term induced by the temporal resolution. In this regime, no matter how large the number of observed trajectories is, the convergence is limited by the too rough (in comparison) temporal resolution.
    \item a ``time-rich'' regime when :
    $$\max \left( \left( \frac{1}{nm} \right)^{\frac{1}{5+b}}, \left( \frac{\ln(n)}{n} \right)^{\frac{1}{b}}, \frac{1}{m} \right)=\left( \frac{\ln(n)}{n} \right)^{\frac{1}{b}}.$$  
    In this case, the convergence rate is \( O\left( \left( \frac{\ln(n)}{n} \right)^{\frac{1}{b}} \right) \), which essentially corresponds to the bias term induced by the sampled initial values. In this regime, no matter how precise the temporal resolution is, the convergence is limited by the too few (in comparison) number of observed trajectories.
    \item a ``balanced'' regime when :
    $$\max \left( \left( \frac{1}{nm} \right)^{\frac{1}{5+b}}, \left( \frac{\ln(n)}{n} \right)^{\frac{1}{b}}, \frac{1}{m} \right)=\left( \frac{1}{nm} \right)^{\frac{1}{5+b}}.$$ 
    In this case, the convergence rate is \( O\left( \left( \frac{1}{nm} \right)^{\frac{1}{5+b}} \right) \), which is optimally obtained through a spatio-temporal bias-variance compromise. In this regime, we benefit from an increase in both the number of observed trajectories and their temporal resolution.
\end{itemize}
Furthermore, under the assumptions we have made, these rates are minimax in the balanced and trajectory-rich regime and, up to logarithmic factors, minimax in the time-rich regime, as stated by the following theorem, whose proof is provided in Section \ref{sec: proof 2}.
\begin{thm}
\label{thm: lowerbound}
Let \( T > 0 \). We have, for $m$ and $n$ sufficiently large :
\[
\inf_{\hat{f}} \sup_{\substack{f \in \operatorname{Lip}(L,M)\\ \mu \in \mathcal{P}(a,b)}} \sup_{x \in \operatorname{Env}_f(\operatorname{supp}(\mu),T)} \mathbb{E} \left[ \|\hat{f}(x) - f(x)\| \right]
\geq K_1 \max \left( \left( \frac{1}{nm} \right)^{\frac{1}{5+b}}, \left( \frac{1}{n} \right)^{\frac{1}{b}}, \frac{1}{m} \right),
\]
where the infimum is taken over all possible estimators of \( f \), and \( K_2 \) is a constant depending only on \( L \), $M$,  \( a \), \( b \), \( D \), \( T \), and \( \sigma \).
\end{thm}
\noindent\textbf{Remark 2.} Several extensions of Theorem \ref{thm: estim-lip-ab} can be considered. First, by examining its proof (see Section \ref{sec:proof 1+crl}), one can observe that the Gaussian assumption on the noise can be relaxed. Indeed, the result remains valid if the family \((\varepsilon_{i,j})_{1\leq i\leq n,1\leq j\leq m}\) is replaced by any family of independent random variables that are sub-Gaussian. Second, if one seeks uniform convergence over the solution envelope, that is, sup-norm convergence, the proof of Theorem~\ref{thm: estim-lip-ab} can be adapted with only minor modifications. In this case, the same rates are recovered, up to additional logarithmic factors (see Appendix~\ref{sec: Exentsion}).\\\\
\textbf{Remark 3.} It is worth noting that, unlike the constrained approaches discussed in the introduction, which require either prior knowledge of the structure supporting the trajectories or its explicit estimation, the construction of our estimator relies on minimal structural information. More precisely, it only requires knowledge of the parameter $b$. Furthermore, our method does not impose strong regularity assumptions on the support of the initial values or trajectories, unlike approaches based on manifold assumptions. For example, \cite{Huang2025} requires the well-definedness of tangent spaces to the underlying manifold.\\\\
\section{Geometric corollaries}\label{sec: geom-coro}
We now derive corollaries under a manifold assumption on the support of the initial values, leveraging the connection between the standardness assumption and the geometry of the support. More precisely, we state that if \(\mu\) has a support \(\mathcal{X}\) that is a well-behaved compact submanifold (without boundary) of \(\mathbb{R}^D\) of dimension \(d\) and admits a density that is lower-bounded by a strictly positive constant on \(\mathcal{X}\), then our estimator achieves convergence rates that depend only on \(d\) and not on \(D\). Although this framework is more restrictive than the one considered in the previous sections, it provides a clearer and more transparent interpretation of our results regarding the role of intrinsic dimension and the curse of dimensionality. 

Before stating this corollary, let us first recall the definition of the reach (of a compact set), a curvature measure introduced by \cite{Fed59}, which has become ubiquitous in manifold learning. Let a compact set \( A \subset \mathbb{R}^D \) and, for all $R\geq 0$, denote $A \oplus R:=\bigcup_{a \in A} B_{2}(a, R)$. The reach of $A$, denoted $\operatorname{reach}(A)$, is the largest $R$ such that each point in $A \oplus R$ has a unique projection onto $A$.
\begin{crl}
\label{crl: manifold assumption}
Let $\mathcal{P}(d,\tau,\varepsilon)$ the collection of probability distributions, $\mu$ , with compact support, $\operatorname{supp}(\mu)$, that are \( d \)-dimensional Riemannian submanifolds of \( \mathbb{R}^D \) with reach \( \tau > 0 \) and admitting a density, $g$, such that $g$ is lower bounded by $\varepsilon$ on $\operatorname{supp}(\mu)$. Choosing $k_{1}$, $k_{2}$, \( h \) and \( k \) as in Theorem \ref{thm: estim-lip-ab} with \( b = d \), we have, for $n$ and $m$ sufficiently large,
\[
\sup_{\substack{f\in \operatorname{Lip}(L,M)\\\mu\in \mathcal{P}(d,\tau,\varepsilon)}}\sup_{x\in \operatorname{Env}_{f}(\operatorname{supp}(\mu),T)}\mathbb{E} \left[ \| \hat{f}(x) - f(x) \| \right] \leq C_2 \max \left( \left( \frac{1}{nm} \right)^{\frac{1}{5+d}}, \left( \frac{\ln(n)}{n} \right)^{\frac{1}{d}}, \frac{1}{m} \right),
\]
where \( C_2 \) is a constant depending only on \( L \), $M$, \( \varepsilon \), \( \tau \), \( d \), \( D \), \( T \), and \( \sigma \).
\end{crl}

 Note that in our context, the assumption of the strict positivity of the reach is not restrictive, as any \( C^2 \)-manifold (in fact, even a \( C^1 \)-manifold with Lipschitz derivatives) has a strictly positive reach. The proof of Corollary \ref{crl: manifold assumption} can be found at the end of Section \ref{sec:proof 1+crl}, where we connect the geometric assumption made here with a standardness assumption. Also note that this rate can also be shown to be minimax, following sensibly the same reasoning that is used for the proof of Theorem \ref{thm: lowerbound} (see section \ref{sec: proof 2}).
\begin{crl}
\label{crl: manifold assumption 2}
Let \( T > 0 \). We have, for $m$ and $n$ sufficiently large :
\[
\inf_{\hat{f}}\sup_{\substack{f\in \operatorname{Lip}(L,M)\\\mu\in \mathcal{P}(d,\tau,\varepsilon)}}\sup_{x\in \operatorname{Env}_{f}(\operatorname{supp}(\mu),T)}\mathbb{E} \left[ \| \hat{f}(x) - f(x) \| \right] \geq K_2 \max \left( \left( \frac{1}{nm} \right)^{\frac{1}{5+d}}, \left( \frac{1}{n} \right)^{\frac{1}{d}}, \frac{1}{m} \right),
\]
where the infimum is taken over all possible estimators of $f$ and \( K_2 \) is a constant depending only on \( L \), $M$, \( \varepsilon \), \( \tau \), \( d \), \( D \), \( T \), and \( \sigma \).
\end{crl}
\textbf{Remark 4.} One may argue that our method relies on knowledge of $d$ (and more generally $b$) to calibrate the parameters $k_1$, $k_2$, $r$, and $k$, which may be unknown in many practical scenarios. Fortunately, the estimation of the intrinsic dimension of a manifold is well-studied, and several efficient estimators have been proposed. For instance, \cite{Kim2016} proposes an estimator that converges at rate $O(\exp(-n\log(n)(1/(D-1)-\varepsilon)))$ for arbitrarily small $\varepsilon>0$, under assumptions very similar to those we make in this section on the distribution of initial conditions. Such an estimator can be used as a preliminary step to estimate the intrinsic dimension from the sampled initial conditions.
\section{Numerical illustrations}
\label{sec: numerical illustration}
\begin{figure}[h!]
    \centering
    \includegraphics[scale=0.38]{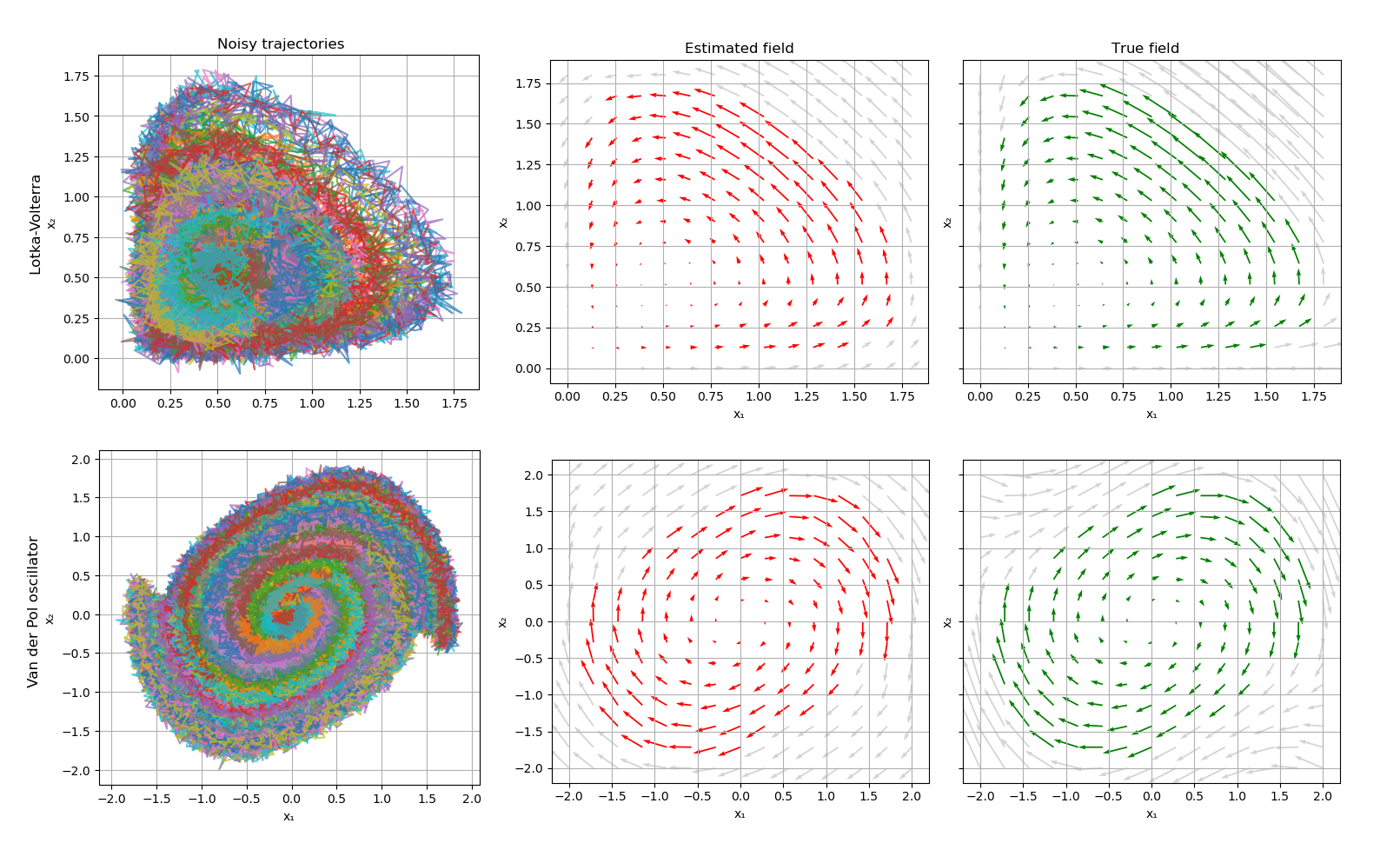}
    \caption{Noisy trajectories along with the estimated and true vector fields for the Van der Pol oscillator and the Lotka--Volterra model. Colored vectors (red and green) indicate positions lying within the respective solution envelopes, $\operatorname{ENV}(\mathcal{X}_1, T_1)$ for the Van der Pol system and $\operatorname{ENV}(\mathcal{X}_2, T_2)$ for the Lotka--Volterra system. For both simulations, we choose $k_1=k=r=10$ and $k_2=7$.}
    \label{fig: simu-1}
\end{figure}

Although the main contributions of this work are theoretical, it is important to assess the practical behavior of our estimator. In this section, we present a series of numerical experiments on synthetic data, which illustrate our theoretical findings while highlighting both the strengths and limitations of our approach.\\\\
\textbf{Qualitative performance on classical 2D models.} We begin by illustrating the qualitative performance of our estimator on two classical (non-linear) two-dimensional dynamical systems: the Van der Pol oscillator and the Lotka-Volterra (predator-prey) model. More precisely, we consider the Van der Pol oscillator defined by:
\[
\begin{cases}
x'_1(t) = x_2(t), \\[6pt]
x'_2(t) = \frac{1}{2}\bigl(1 - x_1(t)^2 \bigr) x_2(t) - x_1(t),
\end{cases}
\]
and the Lotka-Volterra model defined by:
\[
\begin{cases}
x'_1(t) = \frac{1}{2}x_1(t) - x_1(t) x_2(t), \\[6pt]
x'_2(t) = x_1(t) x_2(t) - \frac{1}{2} x_2(t).
\end{cases}
\]
The initial conditions are sampled uniformly along the segment $\mathcal{X}_1=\{(x,x) \mid x \in [-1,1]\}$ for the Van der Pol oscillator and $\mathcal{X}_2=\{(x,x) \mid x \in [1/8,1]\}$ for the Lotka-Volterra model. For both examples, we simulate $n=300$ trajectories, each observed at $m=300$ time points, with additive Gaussian noise of standard deviation $\sigma=0.05$. The time horizon is set to $T_1=4$ for the Van der Pol oscillator and $T_2=10$ for the Lotka-Volterra model. The resulting vector field estimates are displayed in Figure~\ref{fig: simu-1}. The dynamics are well captured within the solution envelopes $\operatorname{ENV}(\mathcal{X}_1, T_1)$ and $\operatorname{ENV}(\mathcal{X}_2, T_2)$, as guaranteed by Theorem~\ref{thm: estim-lip-ab}. Outside these envelopes, the estimation performs poorly; this is particularly evident in the upper-left and lower-left corners of the Van der Pol vector field, and in the upper-left corner of the Lotka-Volterra field. This behavior is expected, as our method relies on locally averaging the flow: regions far from the envelope (and thus far from observed trajectories) suffer from substantial bias. Consequently, while our estimator provides accurate reconstruction within the observed regions, its predictions outside these domains are unreliable.\\\\

\begin{figure}[h!]
  \centering
  \subfloat{\includegraphics[scale=0.38]{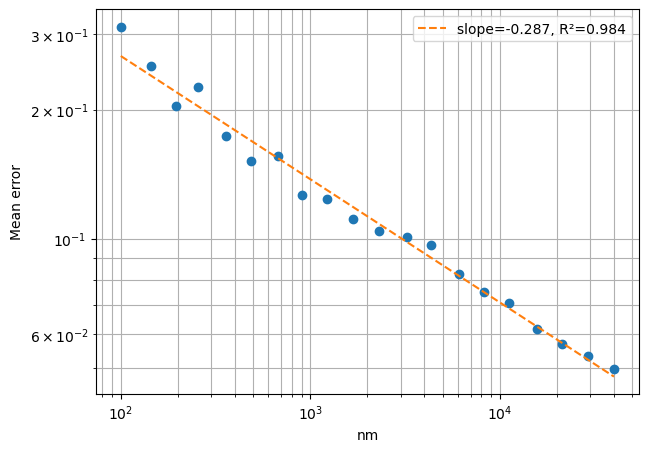}}

  \subfloat{\includegraphics[scale=0.38]{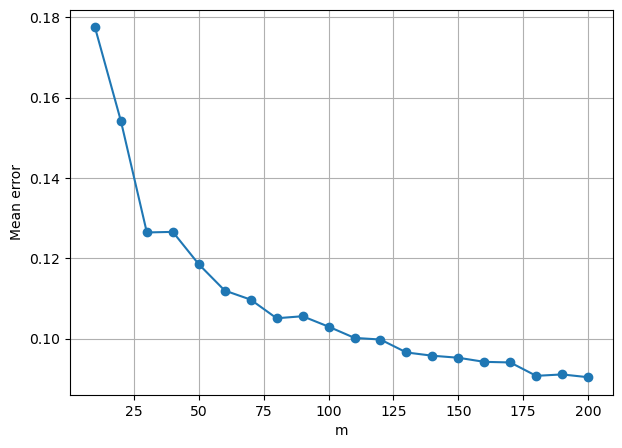}}\hspace{1em}
  \subfloat{\includegraphics[scale=0.38]{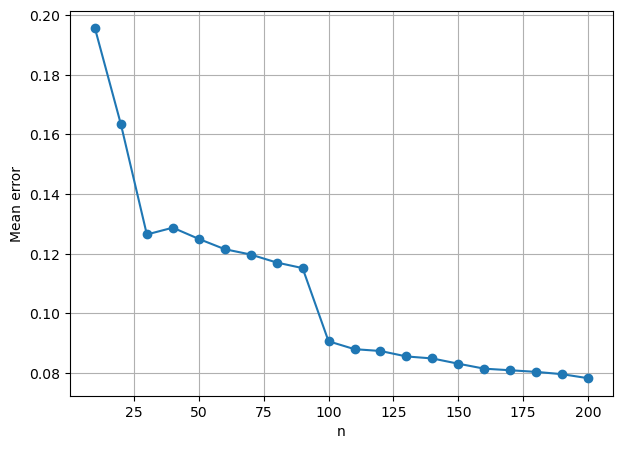}}
  \caption{Mean error as a function of $mn$ in scenario (a) (log--log scale, with the dotted line indicating the linear trend), as a function of $m$ in scenario (b) (linear scale), and as a function of $n$ in scenario (c) (linear scale). Parameters $k_1$, $k_2$, $r$ and $k$ are chosen according to Theorem \ref{thm: estim-lip-ab}.}
  \label{fig: simu-2}
\end{figure}

\noindent\textbf{Different regimes.} We now illustrate the three distinct regimes appearing in Theorems~\ref{thm: estim-lip-ab} and \ref{thm: lowerbound}. To this end, we consider the $D$-dimensional flow defined, for all $1 \leq i \leq D$, by
\begin{equation}\label{eq: Vdp-highdim}
   x'_i(t) = 
\begin{cases} 
x_{i+1}(t), & \text{if $i$ is odd},\\[2mm]
\frac{1}{2} \bigl(1 - x_i(t)^2 \bigr) x_{i+1}(t) - x_i(t), & \text{if $i$ is even}.
\end{cases}
\end{equation}
This system can be seen as a higher-dimensional generalization of the Van der Pol oscillator, where each consecutive pair of dimensions forms a 2-dimensional Van der Pol subsystem. For this simulation, we set $D=6$, $T=2$, and $\sigma=0.05$. Initial conditions are sampled uniformly along the segment \(\mathcal{X} = \{(x,x,\dots,x) \mid x \in [1,2]\}\). We then track the evolution of the mean normalized vector field estimation error $\sum_{i}||\hat{f}(x_i) - f(x_i)||_{\infty} / ||f(x_i)||_{\infty}$ on a $100\times 100$ random grid on the envelope (averaged over $q=10$ repetitions) as a function of $n$ and $m$ in the following three settings:
\begin{itemize}
    \item[(a)] $n = m$, varying from $10$ to $200$,
    \item[(b)] $n = 30$, with $m$ varying from $10$ to $200$,
    \item[(c)] $m = 30$, with $n$ varying from $10$ to $200$.
\end{itemize}
These three scenarios allow us to visualize the three distinct regimes predicted by the theorems. The results, shown in Figure~\ref{fig: simu-2}, confirm the patterns predicted by our theory. When $n = m$, the estimator appears to converge, exhibiting an empirical rate of approximately $(nm)^{-0.28}$, which is faster than the theoretical rate $(nm)^{-1/6}$ guaranteed by Theorem~\ref{thm: estim-lip-ab}. This does not contradict Theorems~\ref{thm: estim-lip-ab} and \ref{thm: lowerbound}, since minimax rates provide a worst-case bound, and faster convergence can be observed in practice, as in this example. In the regime where $n = 30$, the error initially decreases rapidly as $m$ grows, but eventually plateaus, indicating that further increases in $m$ offer little to no improvement. Similarly, when $m = 30$, increasing $n$ reduces the error up to a point, after which additional trajectories provide little to no benefit.\\\\

\textbf{Performance in high ambient dimension.} One of the main motivations behind this work is to study vector field estimation in situations where the initial conditions lie on a low-dimensional structure, while the trajectories evolve in a potentially high-dimensional ambient space. A key property established in Theorem~\ref{thm: estim-lip-ab} and Corollary~\ref{crl: manifold assumption} is that the convergence rates do not depend on the ambient dimension. Another important feature of our estimator is that, since it essentially relies on $k$-nearest neighbors regression techniques, its computational complexity grows only linearly with the ambient dimension. 

\begin{figure}[h]
\centering
\begin{subfigure}{.5\textwidth}
  \centering
  \includegraphics[width=1\linewidth]{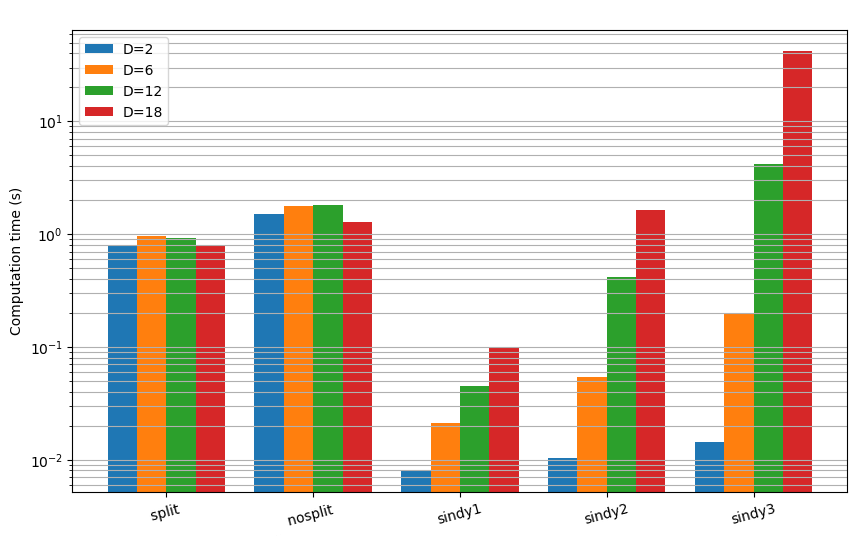}
  \caption{Mean computation time (log scale)}
  \label{fig:sub1}
\end{subfigure}%
\begin{subfigure}{.5\textwidth}
  \centering
  \includegraphics[width=1\linewidth]{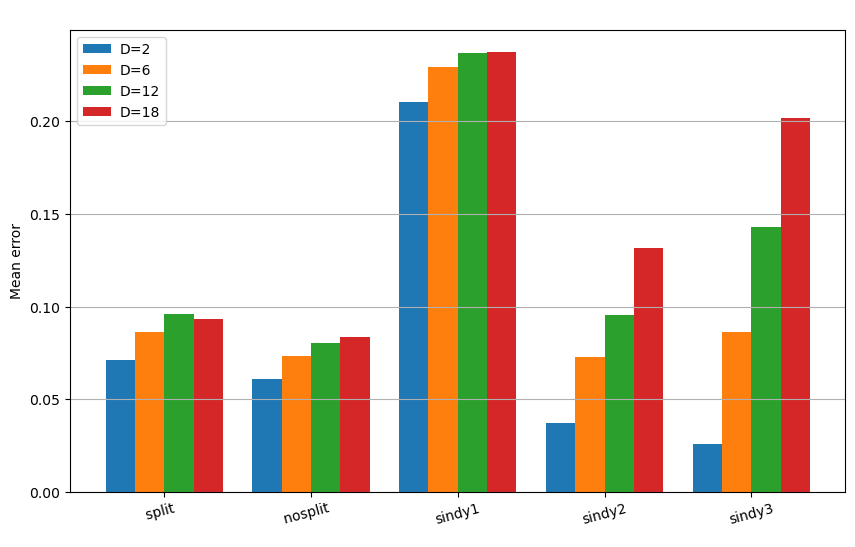}
  \caption{Mean error (linear scale)}
  \label{fig:sub2}
\end{subfigure}
\caption{Mean error and computation time versus ambient dimensions $D=2, 6, 12, 18$ for SINDy estimators (degrees 1, 2, and 3) and for our method (with and without splitting).  Our parameters are chosen according to Theorem \ref{thm: estim-lip-ab}. The thresholds for Sindy estimators are chosen, for each simulation, by minimizing the empirical error within the grid $[0.02, 0.6, 0.1, 0.14, 0.18, 0.22, 0.26, 0.3]$.}
\label{fig: simu-3}
\end{figure}

To illustrate both properties, we generate $n=50$ noisy trajectories ($\sigma = 0.05$) from system~\eqref{eq: Vdp-highdim} with time horizon $T=2$ and $m=100$ time steps, while varying the ambient dimension $D$. We compute both the mean normalized error on the vector field and the mean computation time required to construct the estimator. For comparison, we perform the same evaluation for several SINDy (Sparse Identification of Nonlinear Dynamical systems) estimators~\cite{Brunton2016}. SINDy operates by first estimating the gradients along observed trajectories (similarly to Step~2 of our estimator), then constructing a large library of candidate nonlinear functions of the state variables (e.g., polynomials), and finally applying sparse regression techniques to select a parsimonious subset that best fits the observed dynamics. It is arguably the most popular off-the-shelf algorithm for vector field estimation, notably thanks to its versatile and user-friendly Python implementations~\cite{Desilva2020,Kaptanoglu2022}. For these reasons, it provides a natural point of comparison to evaluate both the statistical accuracy and the computational efficiency of our approach. For the sparse regression component, we adopt the sequential threshold least squares (STLSQ) estimator introduced in \citet{Brunton2016}, which is widely recognized for its robustness to noise. The method proceeds iteratively by performing least squares regression and setting to zero, at each step, coefficients falling below a prescribed threshold. The choice of this threshold is critical yet challenging. In practice, it is often determined via cross-validation over a grid of candidate values. In our simulations, we adopt a simpler strategy: at each iteration, we select the threshold within a predefined grid that minimizes the empirical normalized error. The results of our simulations are reported in Figure~\ref{fig: simu-3}.

First, the results show that the mean normalized error remains relatively stable as the ambient dimension increases, which confirms our theoretical findings. In contrast, SINDy estimators based on degree-two or higher-order polynomial libraries deteriorate significantly with increasing dimension, both statistically and computationally. A similar trend is observed from a computational perspective, as expected: our estimation procedure scales linearly with $D$, whereas the computational cost of SINDy estimators (based on degree-two or higher-order polynomial libraries) grows exponentially with $D$.

Second, although SINDy estimators, particularly those based on degree-two polynomial libraries, perform better in low dimensions, both statistically and computationally, our approach achieves substantially lower mean error and computational cost in higher ambient dimensions. This makes it a compelling alternative to SINDy in such contexts.

Finally, we also tested whether the data-splitting step required in our theoretical analysis is necessary in practice. In this example, omitting the split does not degrade the mean error (and even reduces it slightly), suggesting that while splitting is essential for the proof, it may not be required in practical implementations

\newpage
\section{Conclusion and perspectives}
In this work, we contribute to the statistical understanding of vector field reconstruction from noisy ODEs by introducing and analyzing a flexible regression-like framework from a minimax perspective. We establish minimax rates for pointwise estimation (and for sup-norm estimation, see Appendix~\ref{sec: Exentsion}) under broad distributional assumptions on the initial conditions and mild regularity conditions on the vector field. Our analysis further highlights the relevance of the approach in high ambient dimensions, showing how it helps mitigate the curse of dimensionality, particularly when the initial conditions lie on a low-dimensional manifold. Numerical experiments also demonstrate that our method is computationally efficient—especially in high ambient dimensions, while requiring little prior knowledge about the structure of the solution envelope. In particular, we show that in high-dimensional settings it provides an interesting and simple alternative to classical SINDy estimators, which suffer considerably from the ambient dimension. Although it may not outperform more sophisticated methods in practice, we believe that its attractive statistical and computational properties make it a solid elementary building block for designing more advanced estimators.

Beyond these contributions, our framework opens several promising directions for future research. A natural extension concerns smoother vector fields, such as $C^k$-fields. In nonparametric regression, smoothness assumptions are well known to yield faster convergence rates for function and derivative estimation \citep[see, e.g.,][]{Stone1980}, often via (local) polynomial approximation. This phenomenon also appears in the models of \cite{schötz2024,schötz2025}, where $C^k$-smoothness improves minimax rates. Their analysis, however, crucially relies on strong coverage assumptions: either initial conditions arranged on a uniform grid (Stubble model) or trajectories that densely explore the ambient space (Snake model). By contrast, our setting allows initial conditions and trajectories to lie on lower-dimensional structures, which complicates the extension of their results. Moreover, the local polynomial methods employed in their work suffer from exponential complexity in the ambient dimension, rendering them impractical in high-dimensional contexts. Addressing both of these challenges would represent a substantial advance and motivates further research.
\newpage
\appendix
\section{Proof of Theorem \ref{thm: estim-lip-ab} and Corollary \ref{crl: manifold assumption}}
\label{sec:proof 1+crl}
Let's prove Theorem \ref{thm: estim-lip-ab}. The idea is simply to bound the error component-wise via the following bias-variance decomposition. We denote $f_{s}$ the $s-$th component of $f$, for all $s\in\{1,...,D\}$ and $S_{r}(x)=(X_{i_1(x)},...,X_{i_{r}(x)}, t_{j_1(x)},...,t_{j_r(x)})$, we have :
\begin{align}
\mathbb{E}\left[|\hat{f}_{s}(x)-f_{s}(x)|\right]&= \mathbb{E}\left[\mathbb{E}\left[|\hat{f}_{s}(x)-f_{s}(x)||S_{r}(x)\right]\right]\nonumber\\
    &\leq  \mathbb{E}\left[\sqrt{\operatorname{Var}\left[f_{s}(x)|S_{r}(x)\right]}\right]+ \mathbb{E}\left[\left|\mathbb{E}\left[\hat{f}_{s}(x)|S_{r}(x)\right]-f_{s}(x)\right|\right]\nonumber\\
    &\leq \sqrt{\mathbb{E}\left[\operatorname{Var}\left[f_{s}(x)|S_{r}(x)\right]\right]}+ \mathbb{E}\left[\left|\mathbb{E}\left[\hat{f}_{s}(x)|S_{r}(x)\right]-f_{s}(x)\right|\right]\label{eq:bias-variance decomp}
\end{align}
To bound the two resulting terms, we start by stating two lemmas, proved in Appendix \ref{sec: proof technical lemmas}. Lemma \ref{lmm: Brabanter} is a slight adaptation of Theorem 2 from \cite{Brabanter13}, which bounds the expectation of the variance and the bias of our derivative estimators conditionally on $S_{r}(x)$. 
\begin{lmm}
\label{lmm: Brabanter}
 We have, for all $l\in\{1,...,r\}$ and $s\in\{1,...,D\}$,
\begin{equation}
\label{eq: Brabanter-bias}
\mathbb{E}\left[\left|\mathbb{E}\left[\widehat{f_{s}(\phi(X_{i_{l}(x)},t_{j_{l}(x)}))}|S_{r}(x)\right]-f_{s}(\phi(X_{i_{l}(x)},t_{j_{l}(x)}))\right|\right]\leq C_3\frac{k}{m}
\end{equation}
and 
\begin{equation}
\label{eq: Brabanter-var}
\mathbb{E}\left[\operatorname{Var}\left[\widehat{f_{s}(\phi(X_{i_{l}(x)},t_{j_{l}(x)}))}|S_{r}(x)\right]\right]\leq C_4\frac{m^{2}}{k^{3}}
\end{equation}
with $C_3$ and $C_4$ two constants depending only on $L$, $M$, $T$ and $\sigma$.
\end{lmm}
Lemma \ref{lmm: regression flow} bounds the regression errors for the flow estimator on $Z\times G_T$, ensuring that the obtained convergence rates are faster than the ones stated for vector field estimation.
\begin{lmm}
\label{lmm: regression flow}
Choosing  $k_1$ and $k_2$ as in Theorem \ref{thm: estim-lip-ab}, we have, for $m$ and $n$ sufficiently large :
$$\mathbb{E}\left[\sup_{z\in Z,t\in G_T}\left\|\hat{\phi}(z,t)-\phi(z,t)\right\|\right]\leq C_5 \max\left(\left(\frac{\ln^{2}(nm)}{nm}\right)^{\frac{1}{3+b}},\left(\frac{\ln(n)}{n}\right)^{\frac{1}{b}},\frac{1}{m}\right)$$
with $C_5$ a constant depending only on $L$, $M$, $a$, $b$, $D$, $T$ and $\sigma$. 
\end{lmm}
Note that for all $l\in\{1,...,r\}$ we can write:
$$\widehat{f_{s}(\phi(X_{i_{l}(x)},t_{j_{l}(x)}))}=H(S_{r}(x),(\varepsilon_{i_{l}(x),j})_{1\leq j\leq m}),$$ 
for some measurable function $H$. By construction $S_{r}(x)$ is independent from $(\varepsilon_{i_{l}(x),j})_{1\leq j\leq m}$, and by assumption, $(\varepsilon_{i_{l}(x),j})_{1\leq j\leq m}$ is independent from $(\varepsilon_{i_{l^{'}}(x),j})_{1\leq j\leq m}$, for all $1\leq l^{'}\ne l\leq r$. Thus, $\widehat{f_{s}(\phi(X_{i_{l}(x)},t_{j_{l}(x)}))}$ and $\widehat{f_{s}(\phi(X_{i_{l^{'}}(x)},t_{j_{l^{'}}(x)}))}$ are independent conditionally to $S_{r}(x)$, for all $1\leq l^{'}\ne l\leq r$. Then, from assertion (\ref{eq: Brabanter-var}) of Lemma \ref{lmm: Brabanter}, it follows :
\begin{equation}
\label{eq: var bound}
\mathbb{E}\left[\operatorname{Var}\left[\hat{f}_{s}(x)|S_{r}(x)\right]\right]= \frac{1}{r^{2}}\sum_{l=1}^{r}\mathbb{E}\left[\operatorname{Var}\left[\widehat{f_{s}(\phi(X_{i_{l}(x)},t_{j_{l}(x)}))}|S_{r}(x)\right]\right]\leq C_4\frac{m^2}{rk^{3}}.
\end{equation}
Now, remark that :
\begin{align}
    &\mathbb{E}\left[\left|\mathbb{E}\left[\hat{f}_{s}(x)|S_{r}(x)\right]-f_{s}(x)\right|\right]\nonumber\\
    &\leq\mathbb{E}\left[\left|\mathbb{E}\left[\frac{1}{r}\sum_{l=1}^{r}\widehat{f_{s}(\phi(X_{i_{l}(x)},t_{j_{l}(x)}))}-f_{s}((\phi(X_{i_{l}(x)},t_{j_{l}(x)}))|S_{r}(x)\right]\right|\right]\label{eq: bias bound 1}\\
    &\quad +\mathbb{E}\left[\left|\frac{1}{r}\sum_{l=1}^{r}f_{s}(\phi(X_{i_{l}(x)},t_{j_{l}(x)}))-f_{s}(x)\right|\right].\label{eq: bias bound 2}
\end{align}
By assertion (\ref{eq: Brabanter-bias}) from Lemma \ref{lmm: Brabanter}, we have :
$$(\ref{eq: bias bound 1})\leq\frac{1}{r}\sum_{l=1}^{r}\mathbb{E}\left[\left|\mathbb{E}[\widehat{f_{s}(\phi(X_{i_{l}(x)},t_{j_{l}(x)}))}|S_{r}(x)]-f_{s}(\phi(X_{i_{l}(x)},t_{j_{l}(x)}))\right|\right]\leq C_3\frac{k}{m}.$$
For the second term, as $f\in \operatorname{Lip}(L,M)$, $f_s$ is $L-$Lipschitz and we have:
$$(\ref{eq: bias bound 2})\leq \frac{L}{r}\sum_{l=1}^{r}\mathbb{E}\left[||\phi(X_{i_{l}(x)},t_{j_{l}(x)})-x||\right].$$
Note that there exist $\tilde{x}\in\mathcal{X}$ and $t\in[0,T]$ such that :
$\phi(\tilde{x},t)=x$.  As $f\in \operatorname{Lip}(L,M)$, for all $z\in \mathcal{X}$ and $t'\in[0,T]$, it follows from Grönwall lemma that :
\begin{equation}
\label{eq: deviation bound ODE}
    \|\phi(\tilde{x},t')-\phi(z,t')\|\leq \exp(Lt')\|\tilde{x}-z\|\leq \exp(LT)\|\tilde{x}-z\|.
\end{equation}  
Now, let denote $\tilde{X}=\{\tilde{X}_{1},...,\tilde{X}_{r}\}$ the (ordered) $r-$nearest neighbors in $X[\lfloor n/2\rfloor,n]$ of $\tilde{x}$ and $\tilde{X}^{(r)}\in \operatorname{argmax}\{\min_{k+1\leq j\leq m-k-1}||\hat{\phi}(\tilde{X}_{l},t_j)-x||,\tilde{X}_{l}\in\tilde{X} \}$. By (\ref{eq: deviation bound ODE}) we have:
$$||\phi(\tilde{X}^{(r)},t)-x||\leq ||\phi(\tilde{X}^{(r)},t)-\phi(\tilde{x},t)||\leq \exp(LT)||\tilde{X}_{r}-\tilde{x}||.$$
Let $j\in\operatorname{argmin}\{|t-t_{j}|, k+1\leq j\leq m-k-1\}$, it follows from the definition of the $(X_{i_{l}(x)},t_{j_{l}(x)})$ that, for all $1\leq l\leq r$,
$$||\hat{\phi}(X_{i_{l}(x)},t_{j_{l}(x)})-x||\leq ||\hat{\phi}(\tilde{X}^{(r)},t_{j})-x||.$$
Additionally, note that by the mean value theorem, we have, for all $t,t'\in\mathbb{R}$ and $z\in\mathbb{R}^{d}$:
\begin{equation}
\label{eq: mean value th}
||\phi(z,t)-\phi(z,t')||\leq \sup_{z'\in\mathbb{R}^{d}}||f(z')||\times|t-t'|\leq M |t-t'|.
\end{equation}
It then follows that:
\begin{align*}
&||\phi(X_{i_{l}(x)},t_{j_{l}(x)})-x||\\
&\leq ||\hat{\phi}(X_{i_{l}(x)},t_{j_{l}(x)})-x||+||\hat{\phi}(X_{i_{l}(x)},t_{j_{l}(x)})-\phi(X_{i_{l}(x)},t_{j_{l}(x)})||\\
&\leq ||\hat{\phi}(\tilde{X}^{(r)},t_{j})-x|| + ||\hat{\phi}(X_{i_{l}(x)},t_{j_{l}(x)})-\phi(X_{i_{l}(x)},t_{j_{l}}(x))||\\
&\leq ||\phi(\tilde{X}^{(r)},t_{j})-x||+||\hat{\phi}(\tilde{X}^{(r)},t_{j})-\phi(\tilde{X}^{(r)},t_{j})||+||\hat{\phi}(X_{i_{l}(x)},t_{j_{l}(x)})-\phi(X_{i_{l}(x)},t_{j_{l}(x)})||\\
&\leq ||\phi(\tilde{X}^{(r)},t)-x||+||\phi(\tilde{X}^{(r)},t)-\phi(\tilde{X}^{(r)},t_{j})||\\
&\quad+||\hat{\phi}(\tilde{X}^{(r)},t_{j})-\phi(\tilde{X}^{(r)},t_{j})||+ ||\hat{\phi}(X_{i_{l}(x)},t_{j_{l}(x)})-\phi(X_{i_{l}(x)},t_{j_{l}(x)})||\\
&\leq  \exp(LT)||\tilde{X}_{r}-\tilde{x}||+M|t_{j}-t|\\
&\quad + ||\hat{\phi}(X_{i_{l}(x)},t_{j_{l}(x)})-\phi(X_{i_{l}(x)},t_{j_{l}(x)})||+||\hat{\phi}(\tilde{X}^{(r)},t_{j})-\phi(\tilde{X}^{(r)},t_{j})||\\
&\leq \exp(LT)||\tilde{X}_{r}-\tilde{x}||+2M\frac{k}{m}\\
&\quad + ||\hat{\phi}(X_{i_{l}(x)},t_{j_{l}(x)})-\phi(X_{i_{l}(x)},t_{j_{l}(x)})||+||\hat{\phi}(\tilde{X}^{(r)},t_{j})-\phi(\tilde{X}^{(r)},t_{j})||
\end{align*}
We then have :
\begin{align*}
(\ref{eq: bias bound 2})\leq& \frac{L}{r}\sum_{l=1}^{r}\exp(LT)\mathbb{E}\left[||\tilde{X}_{r}-\tilde{x}||\right]+2M\frac{k}{m}\\
&+\frac{L}{r}\sum_{l=1}^{r}\mathbb{E}\left[||\hat{\phi}(X_{i_{l}(x)},t_{j_{l}(x)})-\phi(X_{i_{l}(x)},t_{j_{l}(x)})||\right]+\mathbb{E}\left[||\hat{\phi}(\tilde{X}^{(r)},t_{j})-\phi(\tilde{X}^{(r)},t_{j})||\right]\\
\leq& L\exp(LT)\mathbb{E}\left[||\tilde{X}_{r}-\tilde{x}||\right]+2LM\frac{k}{m}+L\mathbb{E}\left[||\hat{\phi}(\tilde{X}^{(r)},t_{j})-\phi(\tilde{X}^{(r)},t_{j})||\right]\\
&+\frac{L}{r}\sum_{l=1}^{r}\mathbb{E}\left[||\hat{\phi}(X_{i_{l}(x)},t_{j_{l}(x)})-\phi(X_{i_{l}(x)},t_{j_{l}(x)})||\right]
\end{align*}
Now, we have, by Lemma \ref{lmm: regression flow} :
\begin{align*}
\mathbb{E}\left[||\hat{\phi}(\tilde{X}^{(r)},t_{j})-\phi(\tilde{X}^{(r)},t_{j})||\right]
&\leq \mathbb{E}\left[\sup_{z\in Z,t\in G_T}\left\|\hat{\phi}(z,t)-\phi(z,t)\right\|\right]\\
&\leq C_5 \max\left(\left(\frac{\ln^{2}(nm)}{nm}\right)^{\frac{1}{3+b}},\left(\frac{\ln(n)}{n}\right)^{\frac{1}{b}},\frac{1}{m}\right)   
\end{align*}
and
\begin{align*}
\frac{1}{r}\sum_{l=1}^{r}\mathbb{E}\left[||\hat{\phi}(X_{i_{l}(x)},t_{j_{l}(x)})-\phi(X_{i_{l}(x)},t_{j_{l}(x)})||\right]
&\leq \mathbb{E}\left[\sup_{z\in Z,t\in G_T}\left\|\hat{\phi}(z,t)-\phi(z,t)\right\|\right]\\
&\leq C_5 \max\left(\left(\frac{\ln^{2}(nm)}{nm}\right)^{\frac{1}{3+b}},\left(\frac{\ln(n)}{n}\right)^{\frac{1}{b}},\frac{1}{m}\right).
\end{align*}
It now suffices to bound $\mathbb{E}\left[||\tilde{X}_{r}-\tilde{x}||\right]$. Remark that
$$\mathbb{P}\left(||\tilde{X}_{r}-\tilde{x}||\geq t\right)\leq \mathbb{P}\left(\left|B_{2}(\tilde{x},t)\cap X[\lfloor n/2\rfloor,n]\right|\leq r\right)$$ 
and $\left|B_{2}(\tilde{x},t)\cap X[\lfloor n/2\rfloor,n]\right|$ is a binomial variable of parameter $\left(n-\lfloor n/2\rfloor-1, \mu\left(B_{2}(\tilde{x},t)\right)\right)$, thus by Chernoff bound and using the $(a,b)-$standard assumption on $\mu$ we have :
\begin{align*}
    &\mathbb{P}\left(||\tilde{X}_{r}-\tilde{x}||\geq t\left(\frac{r}{(n-\lfloor n/2\rfloor-1)a}\right)^{\frac{1}{b}}\right)\\
    &\leq \exp\left(-\left(1-\frac{r}{(n-\lfloor n/2\rfloor-1)\mu\left(B_{2}\left(\tilde{x},t\left(\frac{r}{(n-\lfloor n/2\rfloor-1) a}\right)^{\frac{1}{b}}\right)\right)}\right)^{2}\right.\\
    &\quad\quad\quad\quad \left. \times\frac{(n-\lfloor n/2\rfloor-1)\mu\left(B_{2}\left(\tilde{x},t\left(\frac{r}{(n-\lfloor n/2\rfloor-1) a}\right)^{\frac{1}{b}}\right)\right)}{3}\right)\\
    &\leq \exp\left(-\left(1-\frac{1}{t^{b}}\right)^{2}\frac{t^{b}}{3}\right)\\
    &\leq\exp\left(\frac{2}{3}\right)\exp\left(-t^{b}\right).
\end{align*}
Thus, 
\begin{align}
    \mathbb{E}\left[||\tilde{X}_{r}-\tilde{x}||\right]&=\left(\frac{r}{(n-\lfloor n/2\rfloor-1) a}\right)^{\frac{1}{b}}\int_{0}^{+\infty}\mathbb{P}\left(||\tilde{X}_{r}-\tilde{x}||\geq t\left(\frac{r}{(n-\lfloor n/2\rfloor-1) a}\right)^{\frac{1}{b}}\right)dt\nonumber\\
    &\leq \left(\frac{r}{(n-\lfloor n/2\rfloor-1) a}\right)^{\frac{1}{b}}\exp\left(\frac{2}{3}\right)\int_{0}^{+\infty}\exp\left(-t^{b}\right)dt\nonumber\\
    &\leq C(b)\left(\frac{r}{n a}\right)^{\frac{1}{b}}\label{eq: expct NN}
\end{align}
with $C(b)$ a constant depending only on $b$. Hence, we finally obtain :
\begin{align*}
    (\ref{eq: bias bound 2})\leq  2ML\frac{k}{m}+L\exp(LT)C(b)\left(\frac{r}{n a}\right)^{\frac{1}{b}}+2LC_5 \max\left(\left(\frac{\ln^{2}(nm)}{nm}\right)^{\frac{1}{3+b}},\left(\frac{\ln(n)}{n}\right)^{\frac{1}{b}},\frac{1}{m}\right).
\end{align*}
Thus, by (\ref{eq:bias-variance decomp}) :
\begin{align}
 \mathbb{E}\left[|\hat{f}_{s}(x)-f_{s}(x)|\right]
    &\leq (\ref{eq: var bound})+(\ref{eq: bias bound 1})+(\ref{eq: bias bound 2})\nonumber\\
    &\leq C_4^{1/2}\frac{m}{r^{\frac{1}{2}}k^{\frac{3}{2}}}+ C_3\frac{k}{m}+2ML\frac{k}{m}+L\exp(LT)C(b)\left(\frac{r}{n a}\right)^{\frac{1}{b}}\label{eq: intermediate}\\
    &+2LC_5 \max\left(\left(\frac{\ln^{2}(nm)}{nm}\right)^{\frac{1}{3+b}},\left(\frac{\ln(n)}{n}\right)^{\frac{1}{b}},\frac{1}{m}\right)\nonumber.
\end{align}
Now, if $\frac{n^{5/(5+b)}}{m^{b/(5+b)}}>1$
and  $\frac{m^{(4+b)/(5+b)}}{n^{1/(5+b)}}>1$ then $r\simeq \frac{n^{5/(5+b)}}{m^{b/(5+b)}}>1$ and $k\simeq \frac{m^{(4+b)/(5+b)}}{n^{1/(5+b)}}>1$, thus :
\begin{align*}
    (\ref{eq: intermediate})&\leq  C_4^{1/2}\frac{m}{r^{\frac{1}{2}}k^{\frac{3}{2}}}+ C_3\frac{k}{m}+2ML\frac{k}{m}+L\exp(LT)C(b)\left(\frac{r}{na}\right)^{\frac{1}{b}}\\
    &\lesssim \left(\frac{1}{nm}\right)^{\frac{1}{b+5}}
\end{align*}
If $\frac{n^{5/(5+b)}}{m^{b/(5+b)}}<1$ then $m\geq n^{5/b}$, $r=1$ and $k\simeq \frac{m^{(4+b)/(5+b)}}{n^{1/(5+b)}}$, thus :
\begin{align*}
    (\ref{eq: intermediate})&\leq  C_4^{1/2}\frac{m}{k^{\frac{3}{2}}}+ C_3\frac{k}{m}+2ML\frac{k}{m}+L\exp(LT)C(b)\left(\frac{1}{na}\right)^{\frac{1}{b}}\\
    &\lesssim \frac{n^{\frac{3}{2(5+b)}}}{m^{\frac{3(4+b)}{2(5+b)}-1}}+\frac{1}{n^{\frac{1}{5+b}}m^{\frac{1}{5+b}}}+\frac{1}{n^{\frac{1}{b}}}\\
    &\lesssim \frac{1}{n^{\frac{1}{b}}}.
\end{align*}
Finally, if $\frac{m^{(4+b)/(5+b)}}{n^{1/(5+b)}}<1$ then $n\geq m^{4+b}$, $k=1$ and $r\simeq \frac{n^{5/(5+b)}}{m^{b/(5+b)}}$, thus :
\begin{align*}
    (\ref{eq: intermediate})&\leq  C_4^{1/2}\frac{m}{r^{\frac{1}{2}}}+ C_3\frac{1}{m}+2LM\frac{1}{
    m}+L\exp(LT)C(b)\left(\frac{r}{na}\right)^{\frac{1}{b}}\\
    &\lesssim \frac{m^{\frac{b}{2(5+b)}+1}}{n^{\frac{5}{2(5+b)}}}+\frac{1}{m}+\frac{1}{n^{\frac{1}{b+5}}m^{\frac{1}{b+5}}}\\
    &\lesssim \frac{1}{m}.
\end{align*}
Hence, in any case, there exists some constant  \( C_6 \) depending only on \( L \), $M$, \( a \), \( b \), $D$, \( T \), and \( \sigma \) such that :
$$ \mathbb{E}\left[|\hat{f}_{s}(x)-f_{s}(x)|\right] \leq C_6\max\left(\left(\frac{1}{nm}\right)^{\frac{1}{5+b}},\left(\frac{\ln(n)}{n}\right)^{\frac{1}{b}},\frac{1}{m}\right)$$
and Theorem \ref{thm: estim-lip-ab} follows, as $||\hat{f}(x)-f(x)||\leq\sum_{s=1}^{D}|\hat{f}_{s}(x)-f_{s}(x)|$.\qed\\\\
Now for Corollary \ref{crl: manifold assumption}, we suppose that $\mathcal{X}$ has a reach lower bounded by $\tau$. It follows then by Lemma 5.3 from \cite{Niyogi2008} \citep[see also Section 4.2 of][]{ChazalGlisseMichel} that, there exists an absolute constant $C_7$ for all $x\in \mathcal{X}$ and $r<\tau$ :
\begin{align*}
\operatorname{vol}\left(B_{2}(x, r) \cap \mathcal{X}\right) & \geqslant C_7\left(1-\frac{r^2}{4 \tau^2}\right)^{d / 2} r^d \\
& \geqslant C_8 r^d
\end{align*}
with $C_8$ a constant depending on $\tau$ and $d$. Thus, denoting $\mu$ the probability distribution of density $g$ we have, for all $x\in \mathcal{X}$ and $r<\tau$,
\begin{align*}
\mu\left(B_{2}(x,r)\right)&= \mu\left(B_{2}(x,r)\cap \mathcal{X} \right)\\
&=\int_{B_{2}(x,r)\cap \mathcal{X}}g(z)dz\\
&\geq \varepsilon \operatorname{Vol}\left(B_{2}(x,r)\cap \mathcal{X}\right)\\
&\geq \varepsilon C_8 r^d.
\end{align*}
Thus, $\mu$ is $(C_9,d)-$standard, for some constant $C_9$ depending only on $\varepsilon$, $\tau$ and $d$. We can then apply Theorem \ref{thm: estim-lip-ab}.\qed
\section{Proof of Theorem \ref{thm: lowerbound} and Corollary \ref{crl: manifold assumption 2}}
\label{sec: proof 2}
We break the proof of Theorem \ref{thm: lowerbound} in three parts, showing that each term, appearing in the maximum, lower bounds the minimax risk. First, let's prove that :
\begin{equation}
\label{eq: inf 1}
\inf_{\hat{f}}\sup_{f\in \operatorname{Lip}(L,M),
\mu\in \mathcal{P}(a,b)}\sup_{x\in \operatorname{Env}_{f}(X,T)}\mathbb{E}\left[||\hat{f}(x)-f(x)||\right]\gtrsim\left(\frac{1}{nm}\right)^{\frac{1}{5+b}}
\end{equation}
Let $\mu$ a distribution in $\mathcal{P}(a,b)$ with support $\mathcal{X}$ such that for $x_{0}=(0,...,0)$, for all $r>0$, $\mu(B_{2}(x_{0},r))=\min(ar^{b},1)$ and for all $x\in \mathcal{X}$, $x^{t}e_{1}\leq 0$, with $e_{1}$ the first vector in the canonical basis of $\mathbb{R}^{D}$. Let $L,M>0$ and $1\geq h>0$, define the vector fields :
$$f_{0}(x)=\frac{1}{2}\min(M,L)e_{1}$$
and 
$$f_{1}(x)=\frac{1}{2}\min(M,L)\left(1+\left(h-\left\|x-\frac{1}{2}\min(M,L)Te_{1}\right\|_{\infty}\right)_{+}\right)e_{1}.$$
Remark that $f_{0}$ and $f_{1}$ are $L-$Lipschitz and bounded by $M$, $x=(1/2)\min(M,L)Te_{1}$ belong both to $\operatorname{Env}_{f_{0}}(\mathcal{X},T)$ and $\operatorname{Env}_{f_{1}}(\mathcal{X},T)$, and $|f_{1}(x)-f_{0}(x)|=h$. Let $\phi_{0}$ the flow associated with $f_{0}$ and $\phi_{1}$ the flow associated with $f_{1}$. Let $\mathbb{P}_{0}$ the joint distribution of the observations under $f=f_{0}$ and $\mathbb{P}_{1}$ the joint distribution of the observations under $f=f_{1}$. Let $\mathbb{P}_{0}|X$ and $\mathbb{P}_{1}|X$ their distribution conditionally on $X$. We now use a classical variant of Le Cam's two-point argument \citep[see e.g. section 2 of][]{TsybakovBook} that ensures that it suffices to check that if $h=o\left(\left(1/(nm)^{1/(5+b)}\right)\right)$ the Kullback-Liebler divergence ($KL$) of $\mathbb{P}_{0}$ and $\mathbb{P}_{0}$ converges to zero. We have :
\begin{align*}
KL\left(\mathbb{P}_{0},\mathbb{P}_{1}\right)&=\mathbb{E}\left[KL\left(\mathbb{P}_{0}|X,\mathbb{P}_{1}|X\right)\right]\\
    &=\mathbb{E}\left[\sum_{i=1}^{n}\sum_{j=1}^{m}KL\left(\mathcal{N}(\phi_{0}(X_{i},t_{j}),\sigma I_{d}),\mathcal{N}(\phi_{1}(X_{i},t_{j}),\sigma I_{d})\right)\right]\\
    &=\mathbb{E}\left[\sum_{i=1}^{n}\sum_{j=1}^{m}\frac{\|\phi_{0}(X_{i},t_{j})-\phi_{1}(X_{i},t_{j})\|_{2}^{2}}{\sigma^{2}}\right]\\
    &\leq \mathbb{E}\left[\sum_{X_{i}\in [-h,h]^{D}\cap X}\sum_{t_{j}\geq T-h}\frac{\|\phi_{0}(X_{i},t_{j})-\phi_{1}(X_{i},t_{j})\|_{2}^{2}}{\sigma^{2}}\right]\\
    &\leq \mathbb{E}\left[\sum_{X_{i}\in [-h,h]^{D}\cap X}\sum_{t_{j}\geq T-h}\frac{(t_j-(T-h))\min(M,L)^{2}h^2}{4\sigma^{2}}\right]\\
    & \leq \mathbb{E}\left[\sum_{X_{i}\in [-h,h]^{D}\cap X}\sum_{j=1}^{\lfloor hm/T\rfloor}\frac{\min(M,L)^{2}h^{2}j^{2}}{4m^{2}\sigma^{2}}\right]\\
    &\leq \frac{\min(M,L)^{2}h^{2}}{4m^{2}\sigma^{2}}\mathbb{E}\left[\left|X\cap [-h,h]^{D}\right|\right]\sum_{j=1}^{\lfloor hm/T\rfloor}j^{2}\\
    &\leq \frac{\min(M,L)^{2}h^{2}}{4m^{2}\sigma^{2}} \lfloor hm/T\rfloor^{3}\mathbb{E}\left[\left|X\cap [-h,h]^{D}\right|\right]\\
    &\leq  \frac{\min(M,L)^{2}h^{5}m}{2T^{3}\sigma^{2}}\mathbb{E}\left[\left|X\cap [-h,h]^{D}\right|\right].
    \end{align*}
Remark that as the $X_{1},...,X_{n}$ are i.i.d, $|X\cap [-h,h]^{D}|$ follows a binomial distribution of parameters $n$ and $\mathbb{P}\left(X_{1}\in|X\cap [-h,h]^{D}|\right)$, we then have, for $h$ sufficiently small, that:
\begin{align*}
KL\left(\mathbb{P}_{0},\mathbb{P}_{1}\right)&\leq \frac{2\min(M,L)^{2}h^{5}m}{T^{3}\sigma^{2}}\mathbb{E}\left[\left|X\cap [-h,h]^{D}\right|\right]\\
    &\leq \frac{\min(M,L)^{2}h^{5}m}{2T^{3}\sigma^{2}}n\mathbb{P}\left(X_{1}\in B_2(x_{0},\sqrt{d}h)\right)\\
    & \leq \frac{\min(M,L)^{2}h^{5}m}{2T^{3}\sigma^{2}}nad^{b/2}h^{b}=\frac{ad^{b/2}\min(M,L)^{2}}{2T^{3}\sigma^{2}}h^{5+b}nm.
\end{align*}
Thus, if $h=o\left(\left(1/(nm)^{1/(5+b)}\right)\right)$ then $KL\left(\mathbb{P}_{0},\mathbb{P}_{1}\right)\underset{nm\rightarrow\infty}{\longrightarrow}0$ and (\ref{eq: inf 1}) follows.\\\\
Now, let's prove that :
\begin{equation}
\label{eq: inf 3}
\inf_{\hat{f}}\sup_{f\in \operatorname{Lip}(L,M),
\mu\in \mathcal{P}(a,b)}\sup_{x\in \operatorname{Env}_{f}(X,T)}\mathbb{E}\left[||\hat{f}(x)-f(x)||\right]\gtrsim\left(\frac{1}{n}\right)^{\frac{1}{b}}.
\end{equation}
To show this, we place ourselves in an even stronger setting where all the trajectories sampled from $\mu$ are observed entirely (i.e. $m=+\infty$) and without noise (i.e. $\sigma=0$). Let $g$, $f_{0}$ and $f_{1}$ as previously defined. Let $\mathbb{Q}_{0}$ the probability distribution of an observed trajectory under $f=f_{0}$ and  $\mathbb{Q}_{1}$ the probability distribution of an observed trajectory under $f=f_{1}$. We use another classical variant of Le Cam's two-point argument that makes use of the Total Variation distance ($\operatorname{TV}$) between probability distributions, which can be stated in our context as the following inequality \citep[from][]{Yu1997}: 
\begin{align*}
\inf_{\hat{f}}\sup_{f\in \operatorname{Lip}(L,M),
\mu\in \mathcal{P}(a,b)}\sup_{x\in \operatorname{Env}_{f}(X,T)}\mathbb{E}\left[||\hat{f}(x)-f(x)||\right]&\geq \frac{||f_{0}(x)-f_{1}(x)||}{4}\left(1-\operatorname{TV}\left(\mathbb{Q}_{0}^{n},\mathbb{Q}_{1}^{n}\right)\right)\\
&\geq\frac{h}{4}\left(1-n\operatorname{TV}\left(\mathbb{Q}_{0},\mathbb{Q}_{1}\right)\right).
\end{align*}
It then suffices to upper bound the total variation between the distribution $\mathbb{Q}_{0}$ and $\mathbb{Q}_{1}$. Note that for all $x\notin [-h,h]^{D}$ and all $t\in[0,T]$, $\phi_{0}(x,t)=\phi_{1}(x,t)$. Thus, $\mathbb{Q}_{0}$ and $\mathbb{Q}_{1}$ agree on any subsets of the trajectories with initial values not in $[-h,h]^{D}$. Note also that any trajectory is determined by its initial value due to the Cauchy-Lipschitz theorem. Hence, we have :
$$\operatorname{TV}\left(\mathbb{Q}_{0},\mathbb{Q}_{1}\right)\leq 2\int_{[-h,h]^{d}}d\mu\leq 2\mu(B_{2}(0,\sqrt{d}h))\leq 2a(\sqrt{d}h)^{b}.$$
Taking $h=1/(\sqrt{d}(4an)^{1/b})$, we then have :
$$\inf_{\hat{f}}\sup_{f\in \operatorname{Lip}(L,M),
\mu\in \mathcal{P}(a,b)}\sup_{x\in \operatorname{Env}_{f}(X,T)}\mathbb{E}\left[||\hat{f}(x)-f(x)||\right]\geq \frac{1}{8\sqrt{d}(4a)^{\frac{1}{b}}}\left(\frac{1}{n}\right)^{\frac{1}{b}}$$
which gives (\ref{eq: inf 3}).\\\\
Finally, let's prove that :
\begin{equation}
\label{eq: inf 2}
\inf_{\hat{f}}\sup_{f\in \operatorname{Lip}(L,M),
\mu\in \mathcal{P}(a,b)}\sup_{x\in \operatorname{Env}_{f}(X,T)}\mathbb{E}\left[||\hat{f}(x)-f(x)||\right]\gtrsim\frac{1}{m}
\end{equation}
It is indeed a direct corollary of Theorem 3.6 from \cite{schötz2024}, which ensures that for all $x_{0}\in \mathbb{R}^{D}$ there exist two vector fields $f_{0}:\mathbb{R}^{D}\rightarrow \mathbb{R}^{D}$ and $f_{1}:\mathbb{R}^{D}\rightarrow \mathbb{R}^{D}$, belonging to $\operatorname{Lip}(L,M)$ and a constant $C>0$ such that for all $x\in\mathbb{R}^{D}$ and $1\leq j\leq m$ their flows are the same, i.e.:
$$\phi_{0}(x,t_{j})=\phi_{1}(x,t_{j})$$
and such that :
$$||f_{0}(x_{0})-f_{1}(x_{0})||_{2}\geq \frac{C}{m}.$$
Hence, even in the noiseless setting, we have an optimal approximation error lower bounded by $C/m$ and (\ref{eq: inf 2}) follows. Combining (\ref{eq: inf 1}), (\ref{eq: inf 2}), and (\ref{eq: inf 3}), we have proved Theorem \ref{thm: lowerbound}.\qed\\\\
For Corollary \ref{crl: manifold assumption 2}, observe that the proof of (\ref{eq: inf 3}) is entirely deterministic and does not rely on the measure \(\mu\) or its support, so it remains valid in this context. Furthermore, the arguments used in the proofs of (\ref{eq: inf 1}) and (\ref{eq: inf 2}) hold when \(\mathcal{X}\) is a \(d\)-dimensional ball of radius $\tau$, and \(\mu\) is the uniform distribution on \(\mathcal{X}\), in which case, we have $\mu\in \mathcal{P}(d,\tau,\varepsilon)$, for some $\varepsilon$. Consequently, the same reasoning establishes the lower bound stated in Corollary \ref{crl: manifold assumption 2}.\qed
\section{Proof of technical lemmas}
\label{sec: proof technical lemmas}
\textbf{Proof of Lemma \ref{lmm: regression flow}.} Let denote, for all $1\leq p\leq k_{1}$, $i(p,z)$ the index (in $X[1,\lfloor n/2\lfloor]$) associated to $X^{p}(z)$ and, for all $1\leq q\leq k_{2}$, $j(q,t)$ the index associated to $t^{q}(t)$. We have :
\begin{align*}
    \left\|\hat{\phi}(z,t)-\phi(z,t)\right\|&\leq \frac{1}{k_{1}k_{2}}\sum_{p=1}^{k_{1}}\sum_{q=1}^{k_{2}}\left\|\phi(X^{p}(z),t^{q}(t))-\phi(z,t)\right\|\\
    &\quad +\left\|\sigma \frac{1}{k_{1}k_{2}}\sum_{p=1}^{k_{1}}\sum_{q=1}^{k_{2}}\varepsilon_{i(p,z),j(q,t)}\right\|
\end{align*}
Using (\ref{eq: deviation bound ODE}) and (\ref{eq: mean value th}), we have :
\begin{align*}
 &\frac{1}{k_{1}k_{2}}\sum_{p=1}^{k_{1}}\sum_{q=1}^{k_{2}} \left\|\phi(X^{p}(z),t^{q}(t))-\phi(z,t)\right\|\\
 &\leq \frac{1}{k_{1}k_{2}}\sum_{p=1}^{k_{1}}\sum_{q=1}^{k_{2}} \left\|\phi(X^{p}(z),t^{q}(t))-\phi(X^{p}(z),t)\right\|+\left\|\phi(z,t)-\phi(X^{p}(z),t)\right\|\\
   &\leq \frac{1}{k_{1}k_{2}}\sum_{p=1}^{k_{1}}\sum_{q=1}^{k_{2}} M|t-t^{q}(t)|+\exp(LT)||z-X^{p}(z)||\\
   &\leq 2M\frac{k_{2}}{m}+\exp(LT)||z-X^{k_{1}}(z)||
\end{align*}
Hence,
\begin{align}
\mathbb{E}\left[\sup_{z\in Z,t\in G_T}\left\|\hat{\phi}(z,t)-\phi(z,t)\right\|\right]&\leq 2M\frac{k_{2}}{m}+\exp(LT)\mathbb{E}\left[\sup_{z\in Z}||z-X^{k_{1}}(z)||\right]\\
&\quad+\mathbb{E}\left[\sup_{z\in Z,t\in G_T}\left\|\sigma \frac{1}{k_{1}k_{2}}\sum_{p=1}^{k_{1}}\sum_{q=1}^{k_{2}}\varepsilon_{i(p,z),j(q,t)}\right\|\right].
\end{align}
Reasoning as for (\ref{eq: expct NN}) (using an additional union bound argument) and using that $Z$ is independant of $X[1,\lfloor n/2\rfloor]$, one can show that there exists $C(a,b)$ a constant depending on $a$ and $b$ such that :
$$\mathbb{E}\left[\sup_{z\in Z}||z-X^{k_{1}}(z)||\right]\leq C(a,b)\left(\frac{\ln(n)k_{1}}{n}\right)^{\frac{1}{b}}.$$
As $\sigma \frac{1}{k_{1}k_{2}}\sum_{p=1}^{k_{1}}\sum_{q=1}^{k_{2}}\varepsilon_{i(p,z),j(q,t)}$ is a centered Gaussian of covariance matrix $\frac{\sigma^{2}}{k_{1}k_{2}}I_{D}$ conditionnaly to $X=\{X_1,...,X_n\}$, we then have :
\begin{align*}
&\mathbb{E}\left[\sup_{z\in Z,t\in G_T}\left\|\sigma \frac{1}{k_{1}k_{2}}\sum_{p=1}^{k_{1}}\sum_{q=1}^{k_{2}}\varepsilon_{i(p,z),j(q,t)}\right\|\right]\\
&=\mathbb{E}\left[\mathbb{E}\left[\sup_{z\in Z,t\in G_T}\left\|\sigma \frac{1}{k_{1}k_{2}}\sum_{p=1}^{k_{1}}\sum_{q=1}^{k_{2}}\varepsilon_{i(p,z),j(q,t)}\right\||X\right]\right]\\
&\leq 4D\sigma\sqrt{\frac{\ln(nm)}{k_{1}k_{2}}}.
\end{align*}
Thus,
$$\mathbb{E}\left[\sup_{z\in Z,t\in G_T}\left\|\hat{\phi}(z,t)-\phi(z,t)\right\|\right]\leq 2M\frac{k_{2}}{m}+\exp(LT)C(a,b)\left(\frac{\ln(n)k_{1}}{n}\right)^{\frac{1}{b}}+ 4D\sigma\sqrt{\frac{\ln(nm)}{k_{1}k_{2}}}$$
Now, if $\frac{(n/\ln(n))^{3/(b+3)}}{(m/\ln(nm))^{b/(b+3)}}>1$ and $\frac{m^{b+2/(b+3)}}{(n/(\ln(n)\ln(nm)))^{1/(b+3)}}>1$ then :
$$k_{1}\simeq \frac{(n/\ln(n))^{3/(b+3)}}{(m/\ln(nm))^{b/(b+3)}}\text{ and }k_{2}\simeq \frac{m^{b+2/(b+3)}}{(n/(\ln(n)\ln(nm)))^{1/(b+3)}}.$$
Thus: 
$$\mathbb{E}\left[\sup_{z\in Z,t\in G_T}\left\|\hat{\phi}(z,t)-\phi(z,t)\right\|\right]\lesssim \left(\frac{\ln(nm)\ln(n)}{mn}\right)^{\frac{1}{b+3}}\lesssim \left(\frac{\ln^{2}(nm)}{mn}\right)^{\frac{1}{b+3}}.$$
If $\frac{(n/\ln(n))^{3/(b+3)}}{(m/\ln(nm))^{b/(b+3)}}<1$ then $k_1=1$, $k_{2}\simeq \frac{m^{b+2/(b+3)}}{(n/(\ln(n)\ln(nm)))^{1/(b+3)}}$ and $m>\ln(nm)\frac{n^{3/b}}{\ln(n)^{3/b}}$, thus :
\begin{align*}
\mathbb{E}\left[\sup_{z\in Z,t\in G_T}\left\|\hat{\phi}(z,t)-\phi(z,t)\right\|\right]&\leq  2M\frac{k_{2}}{m}+\exp(LT)C(a,b)\left(\frac{\ln(n)}{n}\right)^{\frac{1}{b}}+ 4D\sigma\sqrt{\frac{\ln(nm)}{k_{2}}}\\
&\lesssim \left(\frac{\ln(nm)\ln(n)}{nm}\right)^{\frac{1}{b+3}}+ \left(\frac{\ln(n)}{n}\right)^{\frac{1}{b}}\\
&\quad +\sqrt{\left(\frac{n}{\ln(n)\ln(nm)}\right)^{\frac{1}{b+3}}\frac{\ln(nm)}{m^{\frac{b+2}{b+3}}}}\\
&\lesssim \left(\frac{\ln(n)}{n}\right)^{\frac{1}{b}}
\end{align*}
Finally, if $\frac{m^{b+2/(b+3)}}{(n/(\ln(n)\ln(nm)))^{1/(b+3)}}<1$ then $k_2=1$, $k_{1}\simeq \frac{(n/\ln(n))^{3/(b+3)}}{(m/\ln(nm))^{b/(b+3)}}$ and $n/\ln(n)\ln(nm)>m^{b+2}$, thus :
\begin{align*}
\mathbb{E}\left[\sup_{z\in Z,t\in G_T}\left\|\hat{\phi}(z,t)-\phi(z,t)\right\|\right]&\leq  2M\frac{1}{m}+\exp(LT)C(a,b)\left(\frac{\ln(n)k_{1}}{n}\right)^{\frac{1}{b}}+ 4D\sigma\sqrt{\frac{\ln(nm)}{k_{1}}}\\
&\lesssim \frac{1}{m}+ \left(\frac{\ln(nm)\ln(n)}{nm}\right)^{\frac{1}{b+3}}+\sqrt{\left(\frac{\ln(n)\ln(nm)}{n}\right)^{\frac{3}{b+3}}m^{\frac{b}{b+3}}}\\
&\lesssim \frac{1}{m}
\end{align*}
and Lemma \ref{lmm: regression flow} is proved.\qed\\\\
\textbf{Proof of Lemma \ref{lmm: Brabanter}.} The proof follows, up to slight relaxations and adaptations, from the proof of Theorem 2 from \cite{Brabanter13}. Let's first prove (\ref{eq: Brabanter-bias}). Note that for all $1\leq l\leq r$, $\varepsilon_{i_{l}(x),j_{l}(x)}$ is independent from $S_{r}(x)$, thus, $\mathbb{E}\left[\varepsilon_{i_{l}(x),j_{l}(x)}|S_{r}(x)\right]=\mathbb{E}\left[\varepsilon_{i_{l}(x),j_{l}(x)}\right]=0$. Hence, using that $f$ (and so $f_s$) is $L-$Lipshitz and (\ref{eq: mean value th}), it follows that : 
\begin{align*}
    &\left|\mathbb{E}\left[\widehat{f_{s}(\phi(X_{i_{l}(x)},t_{j_{l}(x)})}|S_{r}(x)\right]-f_{s}(\phi(X_{i_{l}(x)},t_{j_{l}(x)})\right|\\
    &=\left|\sum_{j=1}^{k}w_{j}\frac{\phi_{s}(X_{i_{l}(x)},t_{j_{l}(x)}+j\frac{T}{m})-\phi_{s}(X_{i_{l}(x)},t_{j_{l}(x)}-j\frac{T}{m})}{2j\frac{T}{m}}-f_{s}(\phi_{s}(X_{i_{l}(x)},t_{j_{l}(x)})\right|\\
    &\leq\sum_{j=1}^{k}w_{j}\frac{\left|\phi_{s}(X_{i_{l}(x)},t_{j_{l}(x)}+j\frac{T}{m})-\phi_{s}(X_{i_{l}(x)},t_{j_{l}(x)})-f_{s}(\phi(X_{i_{l}(x)},t_{j_{l}(x)}))j\frac{T}{m}\right|}{2j\frac{T}{m}}\\
    &\quad+\sum_{j=1}^{k}w_{j}\frac{\left|\phi_{s}(X_{i_{l}(x)},t_{j_{l}(x)}-j\frac{T}{m})-\phi_{s}(X_{i_{l}(x)},t_{j_{l}(x)})+f_{s}(\phi(X_{i_{l}(x)},t_{j_{l}(x)}))j\frac{T}{m}\right|}{2j\frac{T}{m}}\\
    &\leq \sum_{j=1}^{k}w_{j}\frac{\int_{0}^{jT/m}|f_s(\phi(X_{i_{l}(x)},t_{j_{l}(x)}+u))-f_{s}(\phi(X_{i_{l}(x)},t_{j_{l}(x)}))|du}{2j\frac{T}{m}}\\
    &\quad+\sum_{j=1}^{k}w_{j}\frac{\int_{0}^{jT/m}|f_s(\phi(X_{i_{l}(x)},t_{j_{l}(x)}))-f_{s}(\phi(X_{i_{l}(x)},t_{j_{l}(x)}-u))|du}{2j\frac{T}{m}}\\
    &\leq L\sum_{j=1}^{k}w_{j}\frac{\int_{0}^{jT/m}||\phi(X_{i_{l}(x)},t_{j_{l}(x)}+u)-\phi(X_{i_{l}(x)},t_{j_{l}(x)})||du}{2j\frac{T}{m}}\\
    &\quad +L\sum_{j=1}^{k}w_{j}\frac{\int_{0}^{jT/m}||\phi(X_{i_{l}(x)},t_{j_{l}(x)}+u)-\phi(X_{i_{l}(x)},t_{j_{l}(x)})||du}{2j\frac{T}{m}}\\
    &\leq LM\sum_{j=1}^{k}w_{j}\left(\frac{\left(j\frac{T}{m}\right)^{2}}{2j\frac{T}{m}}+\frac{\left(j\frac{T}{m}\right)^{2}}{2j\frac{T}{m}}\right)\\
    &\leq LM\sum_{j=1}^{k}w_{j}j\frac{T}{m}\\
    &= \frac{ LMT}{m}\sum_{j=1}^{k}jw_{j}\\
    &=\frac{ LMT}{m}\frac{6}{k(k+1)(2k+1)}\sum_{j=1}^{k}j^{3}\\
    &= \frac{ LMT}{m}\frac{3k^{2}(k+1)^{2}}{2k(k+1)(2k+1)}\\
    &\leq \frac{3}{2} LMT\frac{k}{m}.
\end{align*}
and this still holds, taking the expectation. Now for (\ref{eq: Brabanter-var}), let $\varepsilon^{s}_{i,j}$ denotes the $s-$th component of $\varepsilon_{i,j}$, we have:
\begin{align*}
&\mathbb{E}\left[\operatorname{Var}\left[[\widehat{f_{s}(\phi(X_{i_{l}(x)},t_{j_{l}(x)})}|S_{r}(x)\right]\right]\\
& =\mathbb{E}\left[\operatorname{Var}\left[\sum_{j=1}^{k}w_{j}\sigma\frac{\varepsilon^s_{i_{l}(x),j_{l}(x)+j}-\varepsilon^s_{i_{l}(x),j_{l}(x)-j}}{2j\frac{T}{m}}|S_{r}(x)\right]\right]\\
&\leq \sigma^{2}\operatorname{Var}\left[\sum_{j=1}^{k}w_{j}\frac{\varepsilon^s_{i_{l}(x),j_{l}(x)+j}-\varepsilon^s_{i_{l}(x),j_{l}(x)-j}}{2j\frac{T}{m}}\right]\\
&=\frac{\sigma^{2}m^{2}}{T^{2}}\sum_{j=1}^{k}\frac{w_{j}^{2}}{j^{2}}\\
&=\frac{36\sigma^{2}m^{2}}{T^{2}}\frac{1}{k^{2}(k+1)^{2}(2k+1)^{2}}\sum_{j=1}^{k}j^{2}\\
&=\frac{6\sigma^{2}m^{2}}{T^{2}}\frac{1}{k(k+1)(2k+1)}\leq \frac{6\sigma^{2}}{T^{2}}\frac{m^{2}}{k^{3}}
\end{align*}
and the proof of Lemma \ref{lmm: Brabanter} is complete.\qed
\section{Extension: convergence in sup norm.}
\label{sec: Exentsion}
To rigorously substantiate our remark on the convergence rates of our procedure in the sup norm, we establish the following results. This result guarantees that, up to a $\ln^{2}(nm)$ factor, the sup norm risk achieves the same rates as the pointwise risk. As a result, our proposed strategy remains minimax optimal (up to logarithmic factors), since the lower bound for the pointwise risk directly serves as a lower bound for the sup norm risk.
\begin{thm}
Choosing $k_1$ and $k_2$ as in Theorem \ref{thm: estim-lip-ab},
$$k\simeq \max(m^{(b+4)/(b+5)}((\ln(nm)\ln(n))/n)^{1/(b+5)},1)$$ and 
$$r\simeq\max((n/\ln(n))^{5/(b+5)}(\ln(mn)/m)^{b/(b+5)},1),$$
we have, for $n$ and $m$ sufficiently large :
$$\sup_{\substack{f \in \operatorname{Lip}(L,M)\\\mu \in \mathcal{P}(a,b)}}\mathbb{E} \left[ \sup_{x\in\operatorname{Env}_f(\mathcal{X}, T)}\|\hat{f}(x) - f(x)\| \right]\leq C_{10} \max \left( \left( \frac{\ln^{2}(nm)}{nm} \right)^{\frac{1}{5+b}}, \left( \frac{\ln(n)}{n} \right)^{\frac{1}{b}}, \frac{1}{m} \right)$$
where \( C_{10} \) is a constant depending only on \( L \), $M$, \( a \), \( b \), \( D \), \( T \), and \( \sigma \).
\end{thm}
The proof follows essentially the same steps as for pointwise convergence. First, remark that, for all $x\in \operatorname{Env}_{f}(\mathcal{X},T)$:
$$\hat{f}(x)=\frac{1}{r}\sum_{l=1}^{r}\mathbb{E}\left[\widehat{f\left(\phi(X_{i_{l}(x)},t_{j_{l}(x)})\right)}|S_{r}\right]+\sigma\sum_{u=1}^{k}w_{u}\frac{\varepsilon_{i_l(x),j_l(x)+u}-\varepsilon_{i_l(x),j_l(x)-u}}{2u\frac{T}{m}}.$$
Thus, 
\begin{align}
    &\mathbb{E}\left[\sup_{x\in\operatorname{Env}_{f}(\mathcal{X},T)}\left\|\hat{f}(x)-f(x)\right\|\right]\nonumber\\
    &\leq \frac{1}{r}\sum_{l=1}^{r}\mathbb{E}\left[\sup_{x\in\operatorname{Env}_{f}(\mathcal{X},T)}\left\|\mathbb{E}\left[\widehat{f\left(\phi(X_{i_{l}(x)},t_{j_{l}(x)})\right)}|S_{r}\right]-f\left(\phi(X_{i_{l}(x)},t_{j_{l}(x)})\right)\right\|\right]\label{eq: proba 1}\\
    &\quad+\frac{1}{r}\sum_{l=1}^{r}\mathbb{E}\left[\sup_{x\in\operatorname{Env}_{f}(\mathcal{X},T)}\left\|f\left(\phi(X_{i_{l}(x)},t_{j_{l}(x)})-f(x)\right)\right\|\right]\label{eq: proba 2}\\
    &\quad+\mathbb{E}\left[\sup_{x\in\operatorname{Env}_{f}(\mathcal{X},T)}\left\|\frac{\sigma}{r}\sum_{l=1}^{r}\sum_{u=1}^{k}w_{u}\frac{\varepsilon_{i_l(x),j_l(x)+u}-\varepsilon_{i_l(x),j_l(x)-u}}{2u\frac{T}{m}}\right\|\right]\label{eq: proba 3}.
\end{align}
Let's now bound those three terms. Looking back at the proof of Lemma \ref{lmm: Brabanter} (proof of (\ref{eq: Brabanter-bias})), it follows immediately that :
\begin{equation}
\label{eq: bound infty 1}
(\ref{eq: proba 1})\leq D C_3\frac{k}{m}.
\end{equation}
Now, from Lemma 3 of \cite{Jiang2019}, the number of different $k-$nearest neighbors sets from $N$ points in $\mathbb{R}^{D}$ is bounded by $DN^{D}$. Let denote $\mathcal{M}(Z\times G_T,r)$ the collection of different $r-$nearest neighbors sets from $Z\times G_T$ in $\mathbb{R}^{D+1}$. Using furthermore the independance between $Z$ and $(\varepsilon_{i,j})_{\lfloor n/2 \rfloor< i\leq n, 1\leq j\leq m}$, we then have :
\begin{align}
(\ref{eq: proba 3})&\leq \mathbb{E}\left[\sup_{((i_1,j_1),...,(i_{k},j_{k}))\in \mathcal{M}(Z\times G_T,r) }\left\|\frac{\sigma}{r}\sum_{l=1}^{r}\sum_{u=1}^{k}w_{u}\frac{\varepsilon_{i_j,j_l+u}-\varepsilon_{i_l,j_l-u}}{2u\frac{T}{m}}\right\|\right]\nonumber\\
&\leq 4D\ln^{1/2}(|\mathcal{M}(Z\times G_T,r)|)C_4^{1/2}\frac{m}{r^{1/2}k^{3/2}}\nonumber\\
&\leq 4DC_4^{1/2}(\ln(D+1)+(D+1)\ln(nm))^{1/2}\frac{m}{r^{1/2}k^{3/2}}\label{eq: bound infty 2}.
\end{align}
For the remaining term, reasoning as for the bound we obtained on (\ref{eq: bias bound 2}) we have :
$$(\ref{eq: proba 2})\leq LD\exp(LT)\mathbb{E}\left[\sup_{\tilde{x}\in \mathcal{X}}||\tilde{X}_{r}(\tilde{x})-\tilde{x}||\right]+2LDM\frac{k}{m}+2LD\mathbb{E}\left[\sup_{z\in Z,t\in G_T}||\hat{\phi}(z,t)-\phi(z,t)||\right].$$
By Lemma \ref{lmm: regression flow} :
$$\mathbb{E}\left[\sup_{z\in Z,t\in G_T}||\hat{\phi}(z,t)-\phi(z,t)||\right]\leq C_5 \max\left(\left(\frac{\ln^2(nm)}{nm}\right)^{\frac{1}{3+b}},\left(\frac{\ln(n)}{n}\right)^{\frac{1}{b}},\frac{1}{m}\right).$$
Denote $\mathcal{M}(Z,r)$ the collection of different r-nearest neighbors sets from $Z$ in $\mathcal{X}\subset \mathbb{R}^{D}$. Applying again Lemma 3 from \cite{Jiang2019} and reasoning as in (\ref{eq: expct NN}), we have : 
$$\mathbb{E}\left[\sup_{\tilde{x}\in \mathcal{X}}||\tilde{X}_{r}(\tilde{x})-\tilde{x}||\right]\leq \left(\frac{r\ln(|\mathcal{M}(Z,r)|)}{n}\right)^{\frac{1}{b}}\leq \left(\frac{r(\ln(D)+D\ln(n))}{n}\right)^{\frac{1}{b}}.$$
It follows that : 
\begin{align}
(\ref{eq: proba 2})&\leq LD\exp(LT) \left(\frac{r(\ln(D)+D\ln(n))}{n}\right)^{\frac{1}{b}}+2LDM\frac{k}{m}\nonumber\\
&\quad+2LDC_5 \max\left(\left(\frac{\ln^2(nm)}{nm}\right)^{\frac{1}{3+b}},\left(\frac{\ln(n)}{n}\right)^{\frac{1}{b}},\frac{1}{m}\right).\label{eq: final bound 3} 
\end{align}
Thus, from (\ref{eq: bound infty 1}), (\ref{eq: bound infty 2}) and (\ref{eq: final bound 3}), we have :
\begin{align*}
\mathbb{E}\left[\sup_{x\in\operatorname{Env}_{f}(\mathcal{X},T)}\left\|\hat{f}(x)-f(x)\right\|\right]&\leq (\ref{eq: proba 1})+(\ref{eq: proba 2})+(\ref{eq: proba 3})\\
&\lesssim \frac{k}{m}+\left(\frac{r\ln(n)}{n}\right)^{\frac{1}{b}}+ \ln^{1/2}(nm)\frac{m}{r^{1/2}k^{3/2}}\\
&\quad + \max\left(\left(\frac{\ln^2(nm)}{nm}\right)^{\frac{1}{3+b}},\left(\frac{\ln(n)}{n}\right)^{\frac{1}{b}},\frac{1}{m}\right)\\
&\lesssim \max\left(\left(\frac{\ln^2(nm)}{nm}\right)^{\frac{1}{5+b}},\left(\frac{\ln(n)}{n}\right)^{\frac{1}{b}},\frac{1}{m}\right)
\end{align*}
for the choice 
$$k\simeq \max(m^{(b+4)/(b+5)}((\ln(nm)\ln(n))/n)^{1/(b+5)},1)$$ and 
$$r\simeq\max((n/\ln(n))^{5/(b+5)}(\ln(mn)/m)^{b/(b+5)},1).$$\qed
\newpage
\bibliographystyle{plainnat}
\bibliography{bibliographie} 

@article{McgoffSurvey15,
author = {Kevin McGoff and Sayan Mukherjee and Natesh Pillai},
title = {{Statistical inference for dynamical systems: A review}},
volume = {9},
journal = {Statistics Surveys},
number = {none},
publisher = {Amer. Statist. Assoc., the Bernoulli Soc., the Inst. Math. Statist., and the Statist. Soc. Canada},
pages = {209 -- 252},
year = {2015}
}

@book{RamsaySurvey17,
author = {Ramsay, James and Hooker, Giles},
year = {2017},
month = {01},
pages = {},
title = {Dynamic Data Analysis}
}

@article{DattnerSurvey21,
author = {Dattner, Itai},
title = {Differential equations in data analysis},
journal = {WIREs Computational Statistics},
volume = {13},
number = {6},
pages = {15-34},
year = {2021}
}

@inproceedings{Heinonen18,
  title={Learning unknown ODE models with Gaussian processes},
  author={Markus Heinonen and Çağatay Yıldız and Henrik Mannerstr{\"o}m and Jukka Intosalmi and Harri L{\"a}hdesm{\"a}ki},
  booktitle={International Conference on Machine Learning},
  year={2018}
}

@inproceedings{Chen18,
author = {Chen, Ricky T. Q. and Rubanova, Yulia and Bettencourt, Jesse and Duvenaud, David},
title = {Neural ordinary differential equations},
year = {2018},
booktitle = {Proceedings of the 32nd International Conference on Neural Information Processing Systems},
pages = {6572–6583},
numpages = {12},
series = {NIPS'18}
}

@article{GOTTWALD2021,
title = {Supervised learning from noisy observations: Combining machine-learning techniques with data assimilation},
journal = {Physica D: Nonlinear Phenomena},
volume = {423},
pages = {132911},
year = {2021},
author = {Georg A. Gottwald and Sebastian Reich}
}

@inproceedings{Bhat20,
title={Estimating Vector Fields from Noisy Time Series},
  author={Harish S. Bhat and Majerle Reeves and Ramin Raziperchikolaei},
  booktitle={54th Asilomar Conference on Signals, Systems, and Computers},
  year={2020},
  pages={599-606}
}

@misc{schötz2025,
      title={Nonparametric Estimation of Ordinary Differential Equations: Snake and Stubble}, 
      author={Christof Schötz},
      year={2025}
}

@misc{schötz2024,
      title={Lower Bounds for Nonparametric Estimation of Ordinary Differential Equations}, 
      author={Christof Schötz and Maximilian Siebel},
      year={2024}
}

@article{LAHOUEL2024,
title = {Learning nonparametric ordinary differential equations from noisy data},
journal = {Journal of Computational Physics},
volume = {507},
pages = {112971},
year = {2024},
issn = {0021-9991},
author = {Kamel Lahouel and Michael Wells and Victor Rielly and Ethan Lew and David Lovitz and Bruno M. Jedynak}
}

@article{Brabanter13,
author = {De Brabanter, Kris and De Brabanter, Jos and De Moor, Bart and Gijbels, Ir\`{e}ne},
title = {Derivative estimation with local polynomial fitting},
year = {2013},
volume = {14},
number = {1},
journal = {Journal of Machine Learning Research},
pages = {281–301}
}

@article{Fed59,
    author = {Herbert Federer},
    title = {Curvature measures},
    journal = {Transaction of the American Mathematical Society},
    year = {1959}
}

@article{Niyogi2008,
  title={Finding the Homology of Submanifolds with High Confidence from Random Samples},
  author={Partha Niyogi and Stephen Smale and Shmuel Weinberger},
  journal={Discrete \& Computational Geometry},
  year={2008},
  volume={39},
  pages={419-441}
}

@article{ChazalGlisseMichel,
author = {Chazal, Frédéric and Glisse, Marc and Labruère Chazal, Catherine and Michel, Bertrand},
year = {2014},
pages = {},
title = {Convergence Rates for Persistence Diagram Estimation in Topological Data Analysis},
volume = {1},
journal = {31st International Conference on Machine Learning, ICML 2014}
}

@book{TsybakovBook,
author = {Tsybakov, Alexandre },
title = {Introduction to Nonparametric Estimation},
year = {2008},
publisher = {Springer Publishing Company, Incorporated},
}

@Inbook{Yu1997,
author="Yu, Bin",
title="Assouad, Fano, and Le Cam",
bookTitle="Festschrift for Lucien Le Cam: Research Papers in Probability and Statistics",
year="1997",
pages="423--435"
}

@article{Cuevas2009,
author = {Cuevas, A.},
journal = {Boletín de Estadística e Investigación Operativa. BEIO},
language = {eng},
number = {2},
pages = {71-85},
title = {Set estimation: Another bridge between statistics and geometry.},
volume = {25},
year = {2009},
}

@article{Cuevas2004,
author = {Cuevas, Antonio and Casal, Alberto},
year = {2004},
month = {06},
pages = {340-354},
title = {On boundary estimation},
volume = {36},
journal = {Advances in Applied Probability}
}

@book{Hirsch2013,
  title={Differential equations, dynamical systems, and an introduction to chaos},
  author={Hirsch, Morris W and Smale, Stephen and Devaney, Robert L},
  year={2013},
  publisher={Academic press}
}

@article{Jiang2019, 
title={Non-Asymptotic Uniform Rates of Consistency for k-NN Regression}, 
volume={33},
number={01}, 
journal={Proceedings of the AAAI Conference on Artificial Intelligence},
author={Jiang, Heinrich},
year={2019},
pages={3999-4006} }

@INPROCEEDINGS{Sneed2021,
  author={Sneed, Terry-Ann and Komaee, Arash},
  booktitle={2021 American Control Conference (ACC)}, 
  title={Nonparametric Reconstruction of Vector Fields From Noisy Observations of Their Flow Curves}, 
  year={2021},
  pages={3969-3974}}

@inproceedings{
Huang2025,
title={Learning vector fields of differential equations on manifolds with geometrically constrained operator-valued kernels},
author={Daning Huang and Hanyang He and John Harlim and Yan Li},
booktitle={The Thirteenth International Conference on Learning Representations},
year={2025}
}

@article{Stone1980,
 ISSN = {00905364, 21688966},
 author = {Charles J. Stone},
 journal = {The Annals of Statistics},
 number = {6},
 pages = {1348--1360},
 publisher = {Institute of Mathematical Statistics},
 title = {Optimal Rates of Convergence for Nonparametric Estimators},
 volume = {8},
 year = {1980}
}

@article{Marzouk2024,
author = {Marzouk, Youssef and Ren, Zhi and Wang, Sven and Zech, Jakob},
title = {Distribution learning via neural differential equations: a nonparametric statistical perspective},
year = {2024},
publisher = {JMLR.org},
volume = {25},
number = {1},
issn = {1532-4435},
journal = {J. Mach. Learn. Res.},
articleno = {232},
numpages = {61}
}

@article{
Brunton2016,
author = {Steven L. Brunton  and Joshua L. Proctor  and J. Nathan Kutz },
title = {Discovering governing equations from data by sparse identification of nonlinear dynamical systems},
journal = {Proceedings of the National Academy of Sciences},
volume = {113},
number = {15},
pages = {3932-3937},
year = {2016}}

@article{RUDY2019,
title = {Deep learning of dynamics and signal-noise decomposition with time-stepping constraints},
journal = {Journal of Computational Physics},
volume = {396},
pages = {483-506},
year = {2019},
issn = {0021-9991},
author = {Samuel H. Rudy and J. {Nathan Kutz} and Steven L. Brunton}
}

@article{Desilva2020,
year = {2020},
publisher = {The Open Journal},
volume = {5},
number = {49},
pages = {2104},
author = {Brian de Silva and Kathleen Champion and Markus Quade and Jean-Christophe Loiseau and J. Kutz and Steven Brunton},
title = {PySINDy: A Python package for the sparse identification of nonlinear dynamical systems from data},
journal = {Journal of Open Source Software}
}

@article{Kaptanoglu2022,
year = {2022},
publisher = {The Open Journal},
volume = {7},
number = {69},
pages = {3994},
author = {Alan A. Kaptanoglu and Brian M. de Silva and Urban Fasel and Kadierdan Kaheman and Andy J. Goldschmidt and Jared Callaham and Charles B. Delahunt and Zachary G. Nicolaou and Kathleen Champion and Jean-Christophe Loiseau and J. Nathan Kutz and Steven L. Brunton},
title = {PySINDy: A comprehensive Python package for robust sparse system identification},
journal = {Journal of Open Source Software}
}

@article{Linot2020,
author = {Linot, Alec and Graham, Michael},
year = {2020},
month = {06},
pages = {},
title = {Deep learning to discover and predict dynamics on an inertial manifold},
volume = {101},
journal = {Physical Review E}
}

@article{Vlachas2022,
author = {Vlachas, Pantelis and Arampatzis, Georgios and Uhler, Caroline and Koumoutsakos, Petros},
year = {2022},
month = {04},
pages = {1-8},
title = {Multiscale simulations of complex systems by learning their effective dynamics},
volume = {4},
journal = {Nature Machine Intelligence}
}

@article{Wu2024,
author = {Wu, Tao and Gao, Xiangyun and An, Feng and Sun, Xiaotian and Su, Zhen and Gupta, Shraddha and Gao, Jianxi and Kurths, Jürgen},
year = {2024},
month = {03},
pages = {},
title = {Predicting multiple observations in complex systems through low-dimensional embeddings},
volume = {15},
journal = {Nature Communications}
}

@article{Regazzoni2024,
author = {Regazzoni, Francesco and Pagani, Stefano and Salvador, Matteo and Dede, Luca and Quarteroni, Alfio},
year = {2024},
month = {02},
pages = {},
title = {Learning the intrinsic dynamics of spatio-temporal processes through Latent Dynamics Networks},
volume = {15},
journal = {Nature Communications}
}

@article{Fukami2021,
author = {Fukami, Kai and Murata, Takaaki and Zhang, Kai and Fukagata, Koji},
year = {2021},
month = {11},
title = {Sparse identification of nonlinear dynamics with low-dimensionalized flow representations},
volume = {926},
journal = {Journal of Fluid Mechanics}
}

@article{DAI2016,
  author  = {Wenlin Dai and Tiejun Tong and Marc G. Genton},
  title   = {Optimal Estimation of Derivatives in Nonparametric Regression},
  journal = {Journal of Machine Learning Research},
  year    = {2016},
  volume  = {17},
  number  = {164},
  pages   = {1--25}
}

@misc{Kim2016,
author = {Kim, Jisu and Rinaldo, Alessandro and Wasserman, Larry},
year = {2016},
month = {05},
title = {Minimax Rates for Estimating the Dimension of a Manifold}
}
\section*{Acknowledgements}
The authors would like to thank Eddie Aamari and Marc Hoffmann for our many helpful discussions.
\end{document}